\documentclass{article}
\usepackage{amsmath}
\usepackage{amssymb}

\usepackage[scr=boondoxo]{mathalpha}
  \usepackage{paralist}
  \usepackage[colorlinks=true]{hyperref}
\hypersetup{urlcolor=blue, citecolor=red}

\usepackage{bbm}

\usepackage{geometry}
\newcommand{\ind}[1]{\mathbbm{1}_{#1}}

\usepackage{listings}

\usepackage{booktabs}
\usepackage{makecell}
\usepackage{needspace}

\usepackage{caption}
\usepackage{newfloat}
\DeclareFloatingEnvironment[fileext=frm,placement={!ht},name=Source Code]{listing}
\usepackage{tikz}
\usetikzlibrary{matrix,chains,positioning,decorations.pathreplacing,arrows,math}
\usetikzlibrary{shapes,arrows}
\usepackage{adjustbox}

\tikzset{
  font={\fontsize{7pt}{8}\selectfont}}
  
\usepackage{todonotes}

\usepackage[capitalise]{cleveref} %

\crefname{enumi}{item}{items}
\crefname{subsec}{Section}{Sections}
\crefname{listing}{Source Code}{Source Codes}
\crefname{equation}{}{}

\usepackage{amsthm}

\newtheorem{theorem}{Theorem}[section]
\newtheorem{corollary}[theorem]{Corollary}

\theoremstyle{definition}
\newtheorem{definition}[theorem]{Definition}

\usepackage{verbatim}

\usepackage{cite}

\usepackage[shortlabels,inline]{enumitem}

\usepackage[utf8]{inputenc}
\usepackage[LGR,T1]{fontenc}

\DeclareMathAlphabet{\mathgtt}{LGR}{cmtt}{m}{n}

\usepackage{nicefrac}

\usepackage{xparse}
\usepackage{etoolbox}

\ExplSyntaxOn

\NewDocumentCommand{\enum}{ O{;} m o }
 {
  \my_enum:nnn { #1 } { #2 } { #3 }
 }

\seq_new:N \l__my_enum_seq
\tl_new:N \l__my_enum_item_tl
\prop_const_from_keyval:Nn \l__verbs
{
show = shows ,
imply = implies , 
demonstrate = demonstrates ,
prove = proves ,
establish = establishes ,
ensure = ensures ,
assure = assures
}
\int_new:N \l__number_of_args

\cs_new_protected:Nn \my_enum:nnn
 {
  \seq_set_split:Nnn \l__my_enum_seq { #1 } { #2 }
  \seq_remove_all:Nn \l__my_enum_seq {}
  \int_set_eq:NN \l__number_of_args { \seq_count:N \l__my_enum_seq }
  \seq_use:Nnnn \l__my_enum_seq { ~and~ } { ,~ } { ,~and~ }
  \IfNoValueTF{#3}{}{
    \space
    \int_compare:nNnTF{ \l__number_of_args } < {2}{ \prop_item:Nn \l__verbs {#3} }{ #3 }
  }
 }

\ExplSyntaxOff

 \ExplSyntaxOn

 \seq_new:N \g_cflist_loaded
 \seq_new:N \g_cflist_pending
 
 \NewDocumentCommand{\cfadd}{ m }
 {
   \seq_if_in:NnF \g_cflist_loaded { #1 } {
     \seq_if_in:NnF \g_cflist_pending { #1 } {
       \seq_gput_right:Nn \g_cflist_pending { #1 }
     }
   }
 }

 \NewDocumentCommand{\cfconsiderloaded}{ m }{
   \seq_gput_right:Nn \g_cflist_loaded {#1}
 }
 
 \NewDocumentCommand{\cfremove}{ m }
 {
   \seq_gremove_all:Nn \g_cflist_pending { #1 }
 }

 \NewDocumentCommand{\cfload}{ o }
 {
   \seq_if_empty:NTF \g_cflist_pending {\unskip} {
     (cf.\ \cref{\seq_use:Nn \g_cflist_pending {,}})\IfValueTF{#1}{#1~}{\unskip}
     \seq_gconcat:NNN \g_cflist_loaded \g_cflist_loaded \g_cflist_pending
     \seq_gclear:N \g_cflist_pending
   }
 }
 
 \NewDocumentCommand{\cfclear} {} {
   \seq_gclear:N \g_cflist_loaded
   \seq_gclear:N \g_cflist_pending
 }
 
 \NewDocumentCommand{\cfout}{ o }
 {
   \seq_if_empty:NTF \g_cflist_pending {\unskip} {
     (cf.\ \cref{\seq_use:Nn \g_cflist_pending {,}})\IfValueTF{#1}{#1~}{\unskip}
     \seq_gclear:N \g_cflist_pending
   }
 }
 
 \NewDocumentCommand{\ifnocf} { m } {
   \seq_if_empty:NT \g_cflist_pending { #1 }
 }
 
 \ExplSyntaxOff
 
\ExplSyntaxOn

\NewDocumentEnvironment {athm} {m m o} {%
\IfNoValueTF{#3}{%
\begin{#1}\label{#2}\global\def\loc{#2}%
}{%
\begin{#1}[#3]\label{#2}\global\def\loc{#2}%
}
}{%
\end{#1}%
}

\NewDocumentEnvironment {adef} {m} {%
\begin{definition}\label{#1}\global\def\loc{#1}%
}{%
\end{definition}%
}

\NewDocumentEnvironment{aproof} {} {%
\begin{proof}[Proof~of~\cref{\loc}]%
\bool_gset_false:N \g_finishproof_bool
}{%
\bool_if:NTF \g_finishproof_bool {}
{\finishproofthus}
\end{proof}%
}

\NewDocumentCommand{\finishproofthus} {} {
  \bool_gset_true:N \g_finishproof_bool 
  The~proof~of~\cref{\loc}~is~thus~complete.}
\NewDocumentCommand{\finishproofthis} {} {
  \bool_gset_true:N \g_finishproof_bool 
  This~completes~the~proof~of~\cref{\loc}.}

\ExplSyntaxOff

\newcommand{\lref}[1]{\cref{\loc.#1}}
\newcommand{\llabel}[1]{\label{\loc.#1}}

\ExplSyntaxOn

\NewDocumentEnvironment{flexmath}{ m o }{
  \str_if_eq:noTF {a} {#1} {
    \equation
    \IfValueT{#2}{\label{eq:\loc.#2}}
    \aligned
  } {
    \catcode`&=9
    \renewcommand{\\}{}
    \str_if_eq:noTF {d} {#1} {
      \equation
      \IfValueT{#2}{\label{eq:\loc.#2}}
    } {
      \math
    }
  }
}{
  \str_if_eq:noTF {i} {#1} {
    \endmath
    \catcode`&=4
  } {
    \str_if_eq:noTF {d} {#1} {
      \endequation
    } {
      \endaligned
      \endequation
    }
  }
}

\ExplSyntaxOff

\ExplSyntaxOn

\bool_new:N \g_noteobserve

\NewDocumentCommand{\setnote}{}{
  \bool_gset_true:N \g_noteobserve
}

\NewDocumentCommand{\setobserve}{}{
  \bool_gset_false:N \g_noteobserve
}

\NewDocumentCommand{\nobs}{ o }{
  \IfValueT{#1}{
    \str_if_eq:noTF {note} {#1} {
      \bool_gset_true:N \g_noteobserve
    } {
      \str_if_eq:noTF {Note} {#1} {
        \bool_gset_true:N \g_noteobserve
      } {
        \bool_gset_false:N \g_noteobserve
      }
    }
  }
  \bool_if:nTF { \g_noteobserve } {
    \bool_gset_false:N \g_noteobserve
    note
  } {
    \bool_gset_true:N \g_noteobserve
    observe
  }
  \IfValueF{#1}{~}
}

\NewDocumentCommand{\Nobs}{ o }{
  \IfValueT{#1}{
    \str_if_eq:noTF {note} {#1} {
      \bool_gset_true:N \g_noteobserve
    } {
      \str_if_eq:noTF {Note} {#1} {
        \bool_gset_true:N \g_noteobserve
      } {
        \bool_gset_false:N \g_noteobserve
      }
    }
  }
  \bool_if:nTF { \g_noteobserve } {
    \bool_gset_false:N \g_noteobserve
    Note
  } {
    \bool_gset_true:N \g_noteobserve
    Observe
  }
  \IfValueF{#1}{~}
}

\bool_new:N \g_hencetherefore

\NewDocumentCommand{\hence}{ o }{
  \IfValueT{#1}{
    \str_if_eq:noTF {hence} {#1} {
      \bool_gset_true:N \g_hencetherefore
    } {
      \str_if_eq:noTF {Hence} {#1} {
        \bool_gset_true:N \g_hencetherefore
      } {
        \bool_gset_false:N \g_hencetherefore
      }
    }
  }
  \bool_if:nTF { \g_hencetherefore } {
    \bool_gset_false:N \g_hencetherefore
    hence
  } {
    \bool_gset_true:N \g_hencetherefore
    therefore
  }
  \IfValueF{#1}{~}
}

\NewDocumentCommand{\Hence}{ o }{
  \IfValueT{#1}{
    \str_if_eq:noTF {hence} {#1} {
      \bool_gset_true:N \g_hencetherefore
    } {
      \str_if_eq:noTF {Hence} {#1} {
        \bool_gset_true:N \g_hencetherefore
      } {
        \bool_gset_false:N \g_hencetherefore
      }
    }
  }
  \bool_if:nTF { \g_hencetherefore } {
    \bool_gset_false:N \g_hencetherefore
    Hence~we~obtain~that
  } {
    \bool_gset_true:N \g_hencetherefore
    Therefore~we~obtain~that
  }
  \IfValueF{#1}{~}
}

\int_new:N \g_furthermore

\NewDocumentCommand{\Moreover}{ o o }{
  \IfValueT{#1}{
    \str_case:nn {#1} {
	  {Next} {\int_gset:Nn {\g_furthermore} {0}}      
      {Furthermore} {\int_gset:Nn {\g_furthermore} {1}}
      {Moreover} {\int_gset:Nn {\g_furthermore} {2}}
      {In~addition} {\int_gset:Nn {\g_furthermore} {3}}
      {note} {\bool_gset_true:N \g_noteobserve}
      {observe} {\bool_gset_false:N \g_noteobserve}
    }
    \IfValueT{#2}{
      \str_case:nn {#2} {
	    {Next} {\int_gset:Nn {\g_furthermore} {0}}        
        {Furthermore} {\int_gset:Nn {\g_furthermore} {1}}
        {Moreover} {\int_gset:Nn {\g_furthermore} {2}}
        {In~addition} {\int_gset:Nn {\g_furthermore} {3}}
        {note} {\bool_gset_true:N \g_noteobserve}
        {observe} {\bool_gset_false:N \g_noteobserve}
      }
    }
  }
  \int_case:nn { \int_mod:nn {\g_furthermore} {4} } {
	{ 0 } { Next~\nobs that}    
    { 1 } { Furthermore,~\nobs that}
    { 2 } { Moreover,~\nobs that}
    { 3 } { In~addition,~\nobs that}
  }
  \int_incr:N \g_furthermore
  \IfValueF{#1}{~}
}

\ExplSyntaxOff

\newcommand{\mc}[1]{\mathcal{#1}}
\newcommand{\N}{\mathbb N}
\newcommand{\R}{\mathbb R}

\newcommand{\PP}{\mathbb P}
\newcommand{\Borel}{\mathcal{B}}

\newcommand{\dd}{\mathrm{d}}

\ExplSyntaxOn

\NewDocumentCommand{\mEE}{ o m }{
  \IfValueTF{#1}{
    \str_case:on {#1} {
      {0}{\mathbb E\br{#2}}
      {1}{\mathbb E\br[\big]{#2}}
      {2}{\mathbb E\mkern-1.1mu\br[\Big]{#2}}
      {3}{\mathbb E\mkern-1.3mu\br[\bigg]{#2}}
      {4}{\mathbb E\mkern-1.5mu\br[\Bigg]{#2}}
    }
  } {
    \mathbb E\br{#2}
  }
}

\ExplSyntaxOff

\usepackage{mleftright}

\usepackage{mathtools}
\DeclarePairedDelimiter{\pr}{(}{)}
\DeclarePairedDelimiter{\br}{[}{]}
\DeclarePairedDelimiter{\cu}{\{}{\}}
\DeclarePairedDelimiter{\abs}{\lvert}{\rvert}
\DeclarePairedDelimiter{\norm}{\lVert}{\rVert}

\newcommand{\bbr}[1]{\br[\big]{#1}}

\makeatletter
\newcommand{\mylabel}[2]{#2\def\@currentlabel{#2}\label{#1}}
\makeatother

\newcommand{\E}{\mathbb{E}}

\newcommand{\F}{\mathcal{F}}
\newcommand{\Aff}{\cfadd{def:affine}\mathcal{A}}
\newcommand{\RealV}{\cfadd{def:nets}\mathcal{N}}
\renewcommand{\P}{\mathbb{P}}
\newcommand{\eps}{\varepsilon}

\newcommand{\Python}{\textsc{Python}}

\usepackage{mleftright}

\renewenvironment{pmatrix}{\mleft(\begin{matrix}}{\end{matrix}\mright)}

\usepackage{subcaption}

\usepackage{marginfix}

\title{An overview on deep learning-based approximation\\methods for partial differential equations}

\author{
Christian Beck$^{1,2}$, 
Martin Hutzenthaler$^{3}$,\\
Arnulf Jentzen$^{4,5}$,
Benno Kuckuck$^{6}$
\bigskip
\\
\small{$^1$ Department of Mathematics, ETH Zurich, Zurich, Switzerland}
\smallskip
\\
\small{$^2$ Applied Mathematics: Institute for Analysis and Numerics,}
\vspace{-0.1cm}\\
\small{Faculty of Mathematics and Computer Science,}
\vspace{-0.1cm}\\
\small{University of M\"unster, M\"unster, Germany}
\vspace{-0.1cm}\\
\small{e-mail: \texttt{christian.beck@uni-muenster.de}}
\smallskip
\\
\small{$^3$ Faculty of Mathematics, University of Duisburg-Essen, Essen, Germany}
\vspace{-0.1cm}\\
\small{e-mail: \texttt{martin.hutzenthaler@uni-due.de}}
\smallskip
\\
\small{$^4$ School of Data Science and Shenzhen Research Institute of Big Data,}
\vspace{-0.1cm}\\
\small{The Chinese University of Hong Kong, Shenzhen, China}
\vspace{-0.1cm}\\
\small{e-mail: \texttt{ajentzen@cuhk.edu.cn}}
\smallskip
\\
\small{$^5$ Applied Mathematics: Institute for Analysis and Numerics,}
\vspace{-0.1cm}\\
\small{Faculty of Mathematics and Computer Science,}
\vspace{-0.1cm}\\
\small{University of M{\"u}nster, M\"unster, Germany}
\vspace{-0.1cm}\\
\small{e-mail: \texttt{ajentzen@uni-muenster.de}}
\smallskip
\\
\small{$^6$ Applied Mathematics: Institute for Analysis and Numerics,}
\vspace{-0.1cm}\\
\small{Faculty of Mathematics and Computer Science,}
\vspace{-0.1cm}\\
\small{University of M{\"u}nster, M\"unster, Germany}
\vspace{-0.1cm}\\
\small{e-mail: \texttt{bkuckuck@uni-muenster.de}}
\smallskip
}

\begin{document}
\maketitle

\begin{abstract}
  It is one of the most challenging problems in applied mathematics to 
  approximatively solve high-dimensional partial differential equations (PDEs). 
  Recently, several deep learning-based approximation algorithms
  for attacking
  this problem have been proposed and tested numerically on a 
  number of examples of high-dimensional PDEs. This has given rise
  to a lively field of research in which
  deep learning-based methods and related 
  Monte Carlo methods 
  are applied to the approximation of high-dimensional PDEs.
  In this article we offer an introduction to 
  this field of research by 
  revisiting selected mathematical results 
  related to deep learning approximation methods for PDEs
  and reviewing the main ideas of their proofs.
  We also provide
  a short overview of the recent literature
  in this area of research.
\end{abstract}

\tableofcontents

\section{Introduction}

\label{sec:intro}

Partial differential equations (PDEs) are ubiquitous in
mathematics as tools for
modelling processes in nature
or in man-made complex systems. 
PDEs appear, e.g.,
as Hamilton--Jacobi--Bellman equations
in optimal control problems for describing the value function
associated to the control problem,
as Zakai or Kushner equations
in nonlinear filtering problems for describing the conditional
probability distribution of the state of
the signal process in
the nonlinear filtering problem, 
in models for the approximative valuation 
of financial products such as 
financial derivative contracts,
and in
the approximative description of the distribution of species
in ecosystems to model biodiversity under changing climate conditions.
The PDEs which appear in the above-named applications
are often nonlinear and high-dimensional
where, e.g., in the case of optimal control problems,
the PDE dimension $d\in\N=\{1,2,3,\dots\}$ corresponds to
the number of agents, particles, or resources
in the optimal control problem, 
where, e.g., 
in the case of the approximative
valuation of financial products, the PDE dimension
$d\in\N$ corresponds to the number of financial assets
(such as stocks, commodities, exchange rates, and interest rates)
in the involved hedging portfolio,
and where, e.g.,
in the case of the approximative description
of the distribution of species in ecosystems, 
the PDE dimension $d\in\N$ corresponds to the number
of characteristic traits of the
species in the ecosystem under consideration.

High-dimensional nonlinear PDEs cannot be solved analytically
in nearly all cases
and it is one of the most challenging issues in applied
mathematics to design and analyze approximation
methods for high-dimensional nonlinear PDEs.
Standard approximation methods for nonlinear PDEs, such as
finite difference approximation methods,
finite element approximation methods,
spectral Galerkin approximation methods,
sparse grid approximation methods,
and standard nested Monte Carlo approximation
methods,
suffer from the so-called curse of dimensionality
in the sense that the number of
computational operations of the
employed approximation scheme grows exponentially
in the PDE dimension $d\in\N$
or in the reciprocal $\nicefrac1\eps$ of the 
prescribed approximation
accuracy $\eps\in(0,\infty)$ and it is a very challenging problem
to design and analyze approximation methods
for nonlinear PDEs which overcome the curse of dimensionality
in the sense that the number of computational operations
of the proposed approximation algorithm grows at most
polynomially in the PDE dimension $d\in\N$ and
in the reciprocal $\nicefrac1\eps$ of the
prescribed approximation accuracy $\eps\in(0,\infty)$.

In the last five years, remarkable progress has been made towards
meeting these challenges
and deep learning-based approximation methods have proved particularly
successful in this regard.
Such deep learning-based approximation methods for PDEs have 
first been proposed in the 1990s in the case of low-dimensional
PDEs, cf., e.g., Dissanayake \& Phan-Thien~\cite{dissanayake1994neural},
Lagaris et al.~\cite{lagaris1998artificial},
and Jianyu et al.~\cite{jianyu2003numerical},
and much later, in 2017, in the context of high-dimensional PDEs
in E et al.~\cite{EHanJentzen2017},
Han et al.~\cite{HanJentzenE2018},
Khoo et al.~\cite{khoo2020solving},
Sirignano \& Spiliopoulos~\cite{Sirignano2018dgm},
Beck et al.~\cite{beck2019machine},
E \& Yu~\cite{E2018deep},
and Fujii et al.~\cite{fujii2019asymptotic}. These methods have since been
significantly extended and further studied, and it is one of the
goals of this article to provide a rough overview of the recent
directions this research has taken.

Beside the challenge of solving
PDEs in high dimensions, another area where artificial neural networks
seem to provide a powerful tool is in the context of low-dimensional 
PDEs where not only one fixed solution is desired but instead
a whole family of solutions, depending on the initial condition
and/or additional parameters, needs to be approximated (parametric PDE approximation problems).
Even though the PDE under consideration is low-dimensional,
in this context, a high number of parameters leads to a
high-dimensional parametric approximation problem, which deep learning
approaches seem to be well-suited for;
cf., e.g.~Khoo et al.~\cite{khoo2020solving}
and Li et al.~\cite{li2020fourier} and the references mentioned therein.

In this article we intend to introduce the reader
to the field of research of deep learning-based
approximation of PDEs by revisiting
selected results and the main ideas
of their proofs as well as providing example implementations
in PyTorch and a rough illustration of the capabilities of these methods
as applied to some high-dimensional PDE problems. We present two approximation
methods for nonlinear PDEs, the so-called
deep Galerkin method and the so-called deep splitting method
in some detail and we provide a brief
outline of the relevant literature in the larger field of
 deep learning-based approximation methods
for PDEs.
We note however, that the overview we provide here is by no means
comprehensive. A much longer
article would be needed to do justice
to all the recent developments in this
area.
Our selection is deliberately biased towards
the methods that the authors are most familiar with.
While this unfortunately means that many at least equally deserving ideas
get rather short shrift, we hope that
our introduction encourages the reader to dive deeper
into this fast-growing field of research.

The remainder of this article is organized as follows:
In \cref{sec:linear,sec:nonlinear} of this article we make the idea
of deep learning-based approximation methods
more concrete, where we first focus in \cref{sec:linear} below
on selected deep learning-based approximation methods for linear PDEs and,
thereafter, we briefly sketch in \cref{subsec:dgm,sec:deepsplitting} below some key
mathematical results connected to selected deep learning-based 
approximation methods for nonlinear PDEs, as well as providing 
example implementations.
In \cref{sec:other} we offer a quick overview of the literature on other
deep-learning based approximation methods for PDEs. 
Beside deep learning-based approximation
methods, there are also several other approaches in the scientific 
literature not based on machine learning methods
to approximatively solve high-di\-men\-sional PDEs.
In \cref{sec:othermethods} below, we give a
short overview of the relevant literature regarding these
approaches.
In \cref{sec:simulations} below, we present some simulations
employing the methods introduced in \cref{sec:nonlinear}, illustrating
how these methods can be used to obtain satisfactory results when applied
to certain very high-dimensional PDE problems.
It should be emphasized that though the literature contains a number of
numerical simulations which suggest that
deep learning-based approximation methods may have the capacity to overcome the
curse of dimensionality in the computation of approximate solutions
to high-dimensional PDEs,
there are, however, to date only partial theoretical results 
regarding the conjecture that such methods do indeed
overcome the curse of dimensionality.
These partial results are briefly reviewed in \cref{sec:theoretical} below.
An appendix contains the source code used to obtain the simulation
results laid out in \cref{sec:simulations}.

\section{Deep learning-based approximation methods for linear PDEs}
\label{sec:linear}

In this and the next section of this article (\cref{sec:linear,sec:nonlinear})
we briefly outline how some deep learning-based approximation algorithms 
for PDEs can be derived. In this section we discuss a particular 
deep learning-based approximation algorithm for linear Kolmogorov PDEs 
and in \cref{sec:nonlinear} below 
we treat two deep learning-based approximation algorithms
for nonlinear PDEs, one of which (the deep splitting method
discussed in \cref{sec:deepsplitting} below) is based on the method
presented in this section. 

Linear Kolmogorov PDEs can actually be solved approximately without the 
curse of dimensionality by means of standard Monte Carlo approximation methods. 
However, numerical simulations indicate that deep learning-based
approximation algorithms might be more efficient than
standard Monte Carlo approximation methods when approximating
the solution of a linear Kolmogorov PDE not just at a fixed space-time point
but on an entire region such as on high-dimensional cubes
(cf., e.g., Beck et al.~\cite[Section~4]{beck2021solving}%
).

In this section we sketch the main ideas of the 
deep learning-based approximation algorithm 
for linear Kolmogorov PDEs in~%
\cite[Section~3]{beck2021solving} 
and our presentation here is strongly inspired by the material
in the above-named reference. 
Roughly speaking, the main idea of the deep learning-based approximation algorithm 
in~%
\cite{beck2021solving} is
\begin{enumerate}[(A)]
\item 
\label{item:i_linear_PDEs}
to reformulate the linear Kolmogorov PDE 
whose solution
we intend to approximate 
as a stochastic optimization problem 
on an infinite dimensional space with the unique solution of the 
resulting infinite dimensional stochastic optimization problem 
being the unique solution of the linear Kolmogorov PDE 
which we intend to approximate, 
\item 
\label{item:ii_linear_PDEs}
to approximate the resulting infinite dimensional stochastic 
optimization problem through suitable finite dimensional 
stochastic optimization problems involving 
deep neural networks (DNNs), 
and 
\item 
\label{item:iii_linear_PDEs}
to approximately compute the minimizer
of the resulting finite dimensional 
stochastic optimization problems by means of 
stochastic gradient descent optimization algorithms. 
\end{enumerate}
We now make the above procedure more concrete 
and to further simplify the presentation
we restrict ourselves to linear heat PDEs.

\subsection{Reformulating linear PDEs as infinite dimensional 
stochastic optimization problems}
\label{sec:reformulate}

In \cref{prop:heat_min} below we, roughly speaking, 
reformulate linear heat PDEs as stochastic optimization problems 
on infinite dimensional function spaces (cf.\ \cref{item:i_linear_PDEs} above). 
\cref{prop:heat_min} can be proved by 
an application of the Feynman--Kac formula 
(cf., e.g., \cite[Section~4.4]{karatzasshreve}%
)
and by employing the fact that the expectation of a random variable 
is the best $ L^2 $-approximation of the random variable 
(cf., e.g.,~%
\cite[Lemma~2.1 and Proposition~2.2]{beck2021solving}).
In the scientific literature
\cref{prop:heat_min} is, e.g., proved as~%
\cite[Corollary~2.4]{beck2021solving}.

\begin{theorem}
\label{prop:heat_min}
Let
$ d \in \N $, 
$ T, \rho \in (0,\infty) $, 
$ \varrho = \sqrt{ 2 \rho T } $, 
$ a \in \R $, 
$ b \in (a,\infty) $, 
let 
$ \varphi \colon \R^d \to \R $ be a function,
let 
$
  u 
  \in C^{1,2}([0,T]\times\R^d,\R)
$ 
be a function with at most polynomially growing partial derivatives which satisfies 
for all $t\in [0,T]$, $x\in\R^d$ that 
$ u(0,x) = \varphi(x) $ and 
\begin{equation}
\label{eq:differentialu}
  \tfrac{ \partial u}{\partial t}(t,x) 
  = 
  \rho \, \Delta_x u(t,x)
\end{equation}
let 
$ (\Omega, \F, \PP ) $ be a probability space, 
let
$ \mathbb{W} \colon \Omega \to \R^d $ 
be a standard normal random variable, 
let 
$ \xi \colon \Omega \to [a,b]^d $ be 
a continuous uniformly distributed random variable, 
and
assume that $ \mathbb{W} $ and $ \xi $ are independent.
Then 
\begin{enumerate}[(i)]
\item \label{it:exuniqueSolutions} 
there exists a unique 
$ U \in C( [a,b]^d , \R ) $ such that 
\begin{equation}
\label{eq:stochastic_optimization}
  \E\br[\big]{ \abs{ \varphi( \varrho \mathbb{W} + \xi ) - U(\xi) }^2 }
  = 
  \!\!\inf_{v\in C([a,b]^d,\R)}\!\!\! \E\br[\big]{ \abs{ \varphi( \varrho \mathbb{W}  + \xi ) - v(\xi) }^2 } ,
\end{equation} 
and
\item 
\label{it:UboundaryCond}
it holds for all $ x \in [a,b]^d $ that $ U(x) = u(T,x) $.
\end{enumerate}
\end{theorem}

Roughly speaking, \cref{prop:heat_min} establishes that the solutions 
of the linear heat PDE in \eqref{eq:differentialu} 
can also be viewed as the solutions 
of the stochastic optimization 
problem in \eqref{eq:stochastic_optimization}. 
The deep learning-based approximation method for linear heat PDEs which 
we outline within this section is then based on the approach to approximately 
solve \eqref{eq:stochastic_optimization} in order to obtain 
approximations for the solutions of the PDE in \eqref{eq:differentialu} 
(cf.\ \cref{item:i_linear_PDEs,item:ii_linear_PDEs,item:iii_linear_PDEs} above).

While \cref{prop:heat_min} above recasts the solutions of the PDE in
\eqref{eq:differentialu} at a particular point in time as the solutions
of a stochastic optimization problem, we can also derive from this
a corollary which shows that
the solutions of the PDE over an entire timespan
are similarly the solutions of a stochastic optimization problem.

\begingroup
\begin{corollary}
  \label{cor:kolmogorovtime}
  Let
  $ d \in \N $, 
  $ T, \rho \in (0,\infty) $, 
  $ \varrho = \sqrt{ 2 \rho } $, 
  $ a \in \R $, 
  $ b \in (a,\infty) $, 
  let 
  $ \varphi \colon \R^d \to \R $ be a function,
  let 
  $
    u 
    \in C^{1,2}([0,T]\times\R^d,\R)
  $ 
  be a function with at most polynomially growing partial derivatives which satisfies 
  for all $t\in [0,T]$, $x\in\R^d$ that 
  $ u(0,x) = \varphi(x) $ and 
  \begin{equation}
    \tfrac{ \partial u}{\partial t}(t,x) 
    = 
    \rho \, \Delta_x u(t,x)
    ,
  \end{equation}
  let 
  $ (\Omega, \F, \PP ) $ be a probability space, 
  let
  $ \mathbb{W} \colon \Omega \to \R^d $ 
  be a standard normal random variable, 
  let $\tau\colon \Omega\to[0,T]$ be a continuous uniformly
  distributed random variable,
  let 
  $ \xi \colon \Omega \to [a,b]^d $ be 
  a continuous uniformly distributed random variable, 
  and
  assume that $ \mathbb{W} $, $ \tau $, and $\xi$ are independent.
  Then 
  \begin{enumerate}[(i)]
  \item \label{it:corheat.1}
  there exists a unique 
  $ U \in C( [0,T]\times [a,b]^d , \R ) $ which satisfies that 
  \begin{equation}
  \begin{split}
    &\E\br[\big]{ \abs{ \varphi( \varrho\sqrt{\tau} \mathbb{W} + \xi ) - U(\tau,\xi) }^2 }
    \\&=\!\!\!\!
    \inf_{v\in C([0,T]\times[a,b]^d,\R)}\!\!\!\! \E\br[\big]{ \abs{ \varphi( \varrho\sqrt{\tau}\mathbb{W}  + \xi ) - v(\tau,\xi) }^2 }
  \end{split}
  \end{equation} 
  and
  \item \label{it:corheat.2}
  it holds for all $t\in[0,T]$, $x\in[a,b]^d$ that $ U(t,x) = u(t,x) $.
  \end{enumerate}
\end{corollary}
\renewcommand{\d}[2]{\delta^{#1}_{#2}}
\newcommand{\dee}{\delta}
\newcommand{\e}[1]{\eps_{#1}}
\newcommand{\ee}{\eps}
\begin{proof}[Proof of \cref{cor:kolmogorovtime}]
  Throughout this proof 
  let
    $F\colon C([0,T]\times[a,b]^d,\R)\to[0,\infty]$
  satisfy for all
    $v\in C([0,T]\times[a,b]^d,\R)$
  that
  \begin{equation}
    F(v)
    =
    \E\br[\big]{ \abs{ \varphi( \varrho\sqrt{\tau}\mathbb{W}  + \xi ) - v(\tau,\xi) }^2 }
    .
  \end{equation}
  Note that
    \cref{prop:heat_min}
  establishes that for all
    $v\in C([0,T]\times[a,b]^d,\R)$,
    $s\in[0,T]$
  it holds that
  \begin{equation}
    \label{eq:corheat.2}
    \E\br[\big]{ \abs{ \varphi( \varrho\sqrt{s}\mathbb{W}  + \xi ) - v(s,\xi) }^2 } 
    \geq
    \E\br[\big]{ \abs{ \varphi( \varrho\sqrt{s}\mathbb{W}  + \xi ) - u(s,\xi) }^2 } 
    .
  \end{equation}
  Moreover, observe that 
    the assumption that
      $\mathbb W$, $\tau$, and $\xi$ are independent,
    the assumption that
      $\tau\colon \Omega\to [0,T]$ is continuous uniformly distributed,
    and Fubini's theorem
  ensure that for all
    $v\in C([0,T]\times[a,b]^d,\R)$
  it holds that
  \begin{equation}
    \label{eq:corheat.1}
    F(v)
    =
    \E\br[\big]{ \abs{ \varphi( \varrho\sqrt{\tau}\mathbb{W}  + \xi ) - v(\tau,\xi) }^2 }
    =
    \int_{[0,T]}\E\br[\big]{ \abs{ \varphi( \varrho\sqrt{s}\mathbb{W}  + \xi ) - v(s,\xi) }^2 } \,\dd s
    .
  \end{equation}
    This
    and \eqref{eq:corheat.2}
  prove that for all
    $v\in C([0,T]\times[a,b]^d,\R)$
  it holds that
  \begin{equation}
    F(v)
      \geq
      \int_{[0,T]}\E\bbr{\abs{\varphi(\varrho\sqrt{ s}\mathbb W+\xi)-u(s,\xi)}}\,\dd s
      .
  \end{equation}
  Combining
    this
  with
    \eqref{eq:corheat.1}
  shows that for all
    $v\in C([0,T]\times[a,b]^d,\R)$
  it holds that
  $
    F(v)
    \geq
    F(u)
    $.
  Hence we obtain that
  \begin{equation}
    \label{eq:corheat.3}
    F(u)
    =\!\!\!
    \inf_{v\in C([0,T]\times[a,b]^d,\R)}\!\!\! F(v)
    .
  \end{equation}
    This
    and \eqref{eq:corheat.1}
  demonstrate that for all
    $U\in C([0,T]\times [a,b]^d,\R)$
    with $F(U) = \inf_{v\in C([0,T]\times[a,b]^d,\R)} F(v)$
  it holds that
  \begin{equation}
    \int_{[0,T]}\E\bbr{\abs{\varphi(\varrho\sqrt{ s}\mathbb W+\xi)-U(s,\xi)}}\,\dd s
    =
    \int_{[0,T]}\E\bbr{\abs{\varphi(\varrho\sqrt{ s}\mathbb W+\xi)-u(s,\xi)}}\,\dd s
    .
  \end{equation}
    Combining
      this
    with
     \eqref{eq:corheat.2}
  proves that for all
    $U\in C([0,T]\times [a,b]^d,\R)$
    with $F(U) = \inf_{v\in C([0,T]\times[a,b]^d,\R)} F(v)$
  there exists
    $A\subseteq[0,T]$
    with $\int_A 1\,\dd x=T$
  such that for all
    $s\in A$
  it holds that
  \begin{equation}
    \E\br[\big]{ \abs{ \varphi( \varrho\sqrt{s}\mathbb{W}  + \xi ) - U(s,\xi) }^2 } 
    =
    \E\br[\big]{ \abs{ \varphi( \varrho\sqrt{s}\mathbb{W}  + \xi ) - u(s,\xi) }^2 } 
    .
  \end{equation}
    \cref{prop:heat_min}
    therefore
  ensures that for all
    $U\in C([0,T]\times [a,b]^d,\R)$
    with $F(U) = \inf_{v\in C([0,T]\times[a,b]^d,\R)} F(v)$
  there exists
    $A\subseteq[0,T]$
    with $\int_A 1\,\dd x=T$
  such that for all
    $s\in A$
  it holds that
    $U(s)=u(s)$.
  The fact that
    $u\in C([0,T]\times[a,b]^d,\R)$
    hence
  implies that for all
    $U\in C([0,T]\times [a,b]^d,\R)$
    with $F(U) = \inf_{v\in C([0,T]\times[a,b]^d,\R)} F(v)$
  it holds that
    $U=u$.
  Combining
    this
  with
    \eqref{eq:corheat.3}
  establishes
    \cref{it:corheat.1,it:corheat.2}.
  The proof of \cref{cor:kolmogorovtime} is thus complete.
\end{proof}
\endgroup

\subsection{Mathematical description of deep neural networks (DNNs)}

\newcommand{\rect}[1]{\cfadd{def:rect}{\bf r}_{#1}}
\begingroup
\newcommand{\s}{s}
\newcommand{\netDim}{\mathfrak d}

In order to be able to give a self-contained presentation of
the approximation algorithm sketched in \cref{sec:reformulate} above,
we provide in this subsection a definition of
artificial neural networks suitable for stating
the results in the following sections.
Roughly speaking, a (fully connected feed-forward) artificial neural network is
just a description of a 
function given by alternating compositions of affine linear
functions and particular nonlinear functions.
More formally, 
for all $L\in\{2,3,\dots\}$, $l_0,l_1,\dots,l_L\in\N$, 
$\Psi_1\in C(\R^{l_1},\R^{l_1}),
\Psi_2\in C(\R^{l_2},\R^{l_2}),\allowbreak
\dots,\allowbreak
\Psi_{L-1}\in C(\R^{l_{L-1}},\R^{l_{L-1}})$
and all
affine linear functions 
$A_1\in C(\R^{l_0},\R^{l_1}),
A_2\in C(\R^{l_1},\R^{l_2}),\allowbreak
\ldots,\allowbreak
A_L\in C(\R^{l_{L-1}},\R^{l_L})$
we have that
\begin{equation}
  \label{eq:nnets}
  \R^{l_0}\ni x\mapsto
  (A_L\circ \Psi_{L-1}\circ A_{L-1}\circ\ldots\circ \Psi_1\circ A_1)(x)
  \in \R^{l_L}
\end{equation}
is a realization of a (fully connected feed-forward) artificial neural network
with architecture $(l_0,l_1,\dots,l_L)$ and 
activation functions $(\Psi_1,\Psi_2,\dots,\Psi_{L-1})$.
The following two concepts, \cref{def:affine,def:nets} below, make these
notions more precise.
\Cref{def:affine} introduces a notation for affine linear functions.
\Cref{def:nets} then uses \cref{def:affine} to introduce a
notion of artificial neural networks suitable for our purposes.

\begin{figure}[h!]
  \centering
	\begin{tikzpicture}[shorten >=1pt,-latex,draw=black!100, node distance=\layersep,auto]
			\def\layersep{1.9cm}
			\def\neuronsep{1.2cm}
			\def\ninput{3}
			\def\nhidden{5}
			\def\noutput{3}
			\tikzstyle{every pin edge}=[<-,shorten <=1pt]
			\tikzstyle{neuron}=[circle,fill=black!25,draw=black!75,minimum size=16pt,inner sep=0pt,thick]
			\tikzstyle{input neuron}=[neuron,draw={rgb:red,1;black,1}, fill={rgb:red,1;white,5}];
			\tikzstyle{hidden neuron}=[neuron, draw={rgb:blue,1;black,1}, fill={rgb:blue,1;white,5}];
			\tikzstyle{output neuron}=[neuron,draw={rgb:green,1;black,1}, fill={rgb:green,1;white,5}];
			\tikzstyle{annot} = [text width=9em, text centered]
			\tikzstyle{annot2} = [text width=4em, text centered]

			\tikzmath{
				\ninputl=\ninput-1;
				\nhiddenl=\nhidden-1;
				\noutputl=\noutput-1;
				\nmax=max(\ninput,\nhidden,\noutput);
			}

			\foreach \inode in {1,...,\ninputl}
				\node[input neuron] (I-\inode) at (0,.5*\ninput*\neuronsep-\inode*\neuronsep+\neuronsep) {$\inode$};
			\node (I-dots) at (0,.5*\ninput*\neuronsep-\ninputl*\neuronsep) {$\vdots$};
			\node[input neuron] (I-\ninput) at (0,.5*\ninput*\neuronsep-\ninput*\neuronsep) {$l_0$};

			\foreach \inode in {1,...,\nhiddenl}
				\node[hidden neuron] (H1-\inode) at (\layersep,.5*\nhidden*\neuronsep-\inode*\neuronsep+\neuronsep) {$\inode$};
			\node (H1-dots) at (\layersep,.5*\nhidden*\neuronsep-\nhiddenl*\neuronsep) {$\vdots$};
			\node[hidden neuron] (H1-\nhidden) at (\layersep,.5*\nhidden*\neuronsep-\nhidden*\neuronsep) {$l_1$};

			\foreach \inode in {1,...,\nhiddenl}
				\node[hidden neuron] (H2-\inode) at (2*\layersep,.5*\nhidden*\neuronsep-\inode*\neuronsep+\neuronsep) {$\inode$};
			\node (H2-dots) at (2*\layersep,.5*\nhidden*\neuronsep-\nhiddenl*\neuronsep) {$\vdots$};
			\node[hidden neuron] (H2-\nhidden) at (2*\layersep,.5*\nhidden*\neuronsep-\nhidden*\neuronsep) {$l_2$};

			\foreach \inode in {1,...,\nhiddenl}
				\node (Hdots-\inode) at (3*\layersep,.5*\nhidden*\neuronsep-\inode*\neuronsep+\neuronsep) {$\cdots$};
			\node (Hdots-dots) at (3*\layersep,.5*\nhidden*\neuronsep-\nhiddenl*\neuronsep) {$\ddots$};
			\node (Hdots-\nhidden) at (3*\layersep,.5*\nhidden*\neuronsep-\nhidden*\neuronsep) {$\cdots$};

			\foreach \inode in {1,...,\nhiddenl}
				\node[hidden neuron] (H3-\inode) at (4*\layersep,.5*\nhidden*\neuronsep-\inode*\neuronsep+\neuronsep) {$\inode$};
			\node (H3-dots) at (4*\layersep,.5*\nhidden*\neuronsep-\nhiddenl*\neuronsep) {$\vdots$};
			\node[hidden neuron] (H3-\nhidden) at (4*\layersep,.5*\nhidden*\neuronsep-\nhidden*\neuronsep) {$l_{L-1}$};

			\foreach \inode in {1,...,\noutputl}
				\node[output neuron] (O-\inode) at (5*\layersep,.5*\noutput*\neuronsep-\inode*\neuronsep+\neuronsep) {$\inode$};
			\node (O-dots) at (5*\layersep,.5*\noutput*\neuronsep-\noutputl*\neuronsep) {$\vdots$};
			\node[output neuron] (O-\noutput) at (5*\layersep,.5*\noutput*\neuronsep-\noutput*\neuronsep) {$l_L$};

			\foreach \inode in {1,...,\ninput}
				\foreach \hnode in {1,...,\nhidden}
					\path (I-\inode) edge [draw=black!30] (H1-\hnode);

			\foreach \inode in {1,...,\nhidden}
				\foreach \hnode in {1,...,\nhidden}
					\path (H1-\inode) edge [draw=black!30] (H2-\hnode);

			\foreach \inode in {1,...,\nhidden}
				\foreach \hnode in {1,...,\nhidden}
					\path (H2-\inode) edge [-,draw=black!30] (Hdots-\hnode);

			\foreach \inode in {1,...,\nhidden}
				\foreach \hnode in {1,...,\nhidden}
					\path (Hdots-\inode) edge [draw=black!30] (H3-\hnode);

			\foreach \inode in {1,...,\nhidden}
				\foreach \hnode in {1,...,\noutput}
					\path (H3-\inode) edge [draw=black!30] (O-\hnode);

			\node[annot] (input) at (0,.5*\nmax*\neuronsep+.5*\neuronsep) {Input layer\\(1st layer)}; 
			\node[annot] (hidden1) at (\layersep,.5*\nmax*\neuronsep+0.85*\neuronsep) {1st hidden layer\\(2nd layer)}; 
			\node[annot] (hidden2) at (2*\layersep,.5*\nmax*\neuronsep+.5*\neuronsep) {2nd hidden layer\\(3rd layer)}; 
			\node[annot] (dots) at (3*\layersep,.5*\nmax*\neuronsep+.5*\neuronsep) {$\cdots$}; 
			\node[annot] (hidden3) at (4*\layersep,.5*\nmax*\neuronsep+.85*\neuronsep) {$(L-1)$-th hidden layer\\($L$-th layer)}; 
			\node[annot] (output) at (5*\layersep,.5*\nmax*\neuronsep+.5*\neuronsep) {Output layer\\($(L+1)$-th layer)}; 
	\end{tikzpicture}
\caption{\label{fig:shallowANN}Graphical illustration of a fully-connected 
feedforward artificial neural network consisting of
$L\in\N$ affine linear transformations (i.e., consisting of $L+1$ layers: one input layer, $L-1$ hidden layers, and one output layer) 
with $l_0\in\N$ neurons on the input layer, with
$l_1\in\N$ neurons on the first hidden layer,
with $l_2\in\N$ neurons on the second hidden layer,
$\dots$, with $l_{L-1}$ neurons on the $(L-1)$-th hidden layer,
and with $l_L$ neurons in the output layer.}
\end{figure}
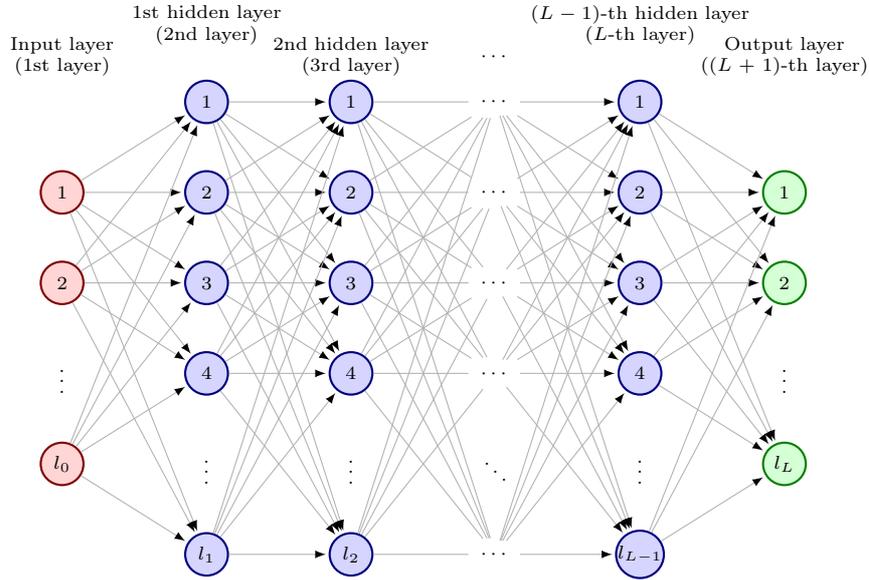

\begin{definition}[Affine linear functions]\label{def:affine}
  Let $\mathfrak d,m,n \in \N$, $s \in \N_0$, 
$ \theta = ( \theta_1, \theta_2, \dots,\allowbreak \theta_\netDim ) \in \R^\netDim $ 
satisfy $\netDim \geq \s + m n + m$.
Then we denote by $\Aff_{m,n}^{\theta, \s}\colon\R^n\to\R^m$ the function which satisfies 
for all $x = (x_1,x_2,\ldots, x_{n}) \in \R^{n}$ that
\begin{equation}
  \label{eq:affine}
\begin{split}
   &\Aff_{m,n}^{\theta,\s}( x ) 
= 
  \begin{pmatrix}
      \theta_{ \s + 1 }
    &
      \theta_{ \s + 2 }
    &
      \cdots
    &
      \theta_{ \s + n }
    \\
      \theta_{ \s + n + 1 }
    &
      \theta_{ \s + n + 2 }
    &
      \cdots
    &
      \theta_{ \s + 2 n }
    \\
      \theta_{ \s + 2 n + 1 }
    &
      \theta_{ \s + 2 n + 2 }
    &
      \cdots
    &
      \theta_{ \s + 3 n }
    \\
      \vdots
    &
      \vdots
    &
      \ddots
    &
      \vdots
    \\
      \theta_{ \s + ( m - 1 ) n + 1 }
    &
      \theta_{ \s + ( m - 1 ) n + 2 }
    &
      \cdots
    &
      \theta_{ \s + m n }
    \end{pmatrix}
    \begin{pmatrix}
      x_1
    \\
      x_2
    \\
      x_3
    \\
      \vdots 
    \\
      x_{n}
    \end{pmatrix}
  +
    \begin{pmatrix}
      \theta_{ \s + m n + 1 }
    \\
      \theta_{ \s + m n + 2 }
    \\
      \theta_{ \s + m n + 3 }
    \\
      \vdots 
    \\
      \theta_{ \s + m n + m }
    \end{pmatrix}
.
\end{split}
\end{equation}
\end{definition}

Observe that for all $\mathfrak d,m,n \in \N$, $s \in \N_0$,
$\theta\in\R^{\mathfrak d}$ the definition
of $\Aff_{m,n}^{\theta, \s}\colon\R^n\to\R^m$ in \cref{def:affine} above
refers to 
the components $\theta_{s+1},\theta_{s+2},\dots,\theta_{s+mn+m}$.
The condition that
$\mathfrak d\geq \s+mn+m$ hence ensures that
the definition of $\Aff_{m,n}^{\theta, \s}\colon\R^n\to\R^m$ in
\eqref{eq:affine} above make sense.
In \cref{def:nets} below, we now use the affine linear functions introduced in
\cref{def:affine} above to define the realizations of artificial neural networks.

\cfclear
\cfconsiderloaded{def:nets}
\begin{definition}[Realizations of artificial neural networks]\label{def:nets}
  Let
$L\in\{2,3,\dots\}$, $l_0,l_1,\dots,l_L
\in \N 
$,
$\Psi_1\in C(\R^{l_1},\R^{l_1}),\allowbreak
\Psi_2\in C(\R^{l_2},\R^{l_2}),\allowbreak
\dots,\allowbreak
\Psi_{L-1}\in C(\R^{l_{L-1}},\R^{l_{L-1}})$
and for every
$j\in\{1,2,\dots,L+1\}$
let
$\mathfrak d_j\in\N_0$
satisfy $\mathfrak d_j=\sum_{ k = 1 }^{ j-1 } l_k (l_{k-1} + 1)$.
Then we denote by
$
  \mathcal N^{l_0,l_1,\dots,l_L}_{\Psi_1,\Psi_2,\dots,\Psi_{L-1}}
  \colon
  \R^{\mathfrak d_{L+1}}\to C(\R^{l_0},\R^{l_L})
$
the function which satisfies for all
$\theta\in \R^{\mathfrak d_{L+1}}$
that
\begin{equation}
  \label{eq:paramFFNN}
  \pr[\big]{\mathcal N^{l_0,l_1,\dots,l_L}_{\Psi_1,\Psi_2,\dots,\Psi_{L-1}}}(\theta)
  =
  \Aff_{l_{L-1},l_L}^{\theta,\mathfrak d_L}\circ \Psi_{L-1}
  \circ\Aff_{l_{L-2},l_{L-1}}^{\theta,\mathfrak d_{L-1}}
  \circ\ldots
  \circ \Psi_1\circ \Aff_{l_0,l_1}^{\theta,\mathfrak d_1}
\end{equation}
\cfload.
\end{definition}
We observe that for all
$L\in\N$, $l_0,l_1,l_2,\dots,l_L\in\N$,
$\Psi_1\in C(\R^{l_1},\R^{l_1}),
\Psi_2\in C(\R^{l_2},\R^{l_2}),\allowbreak
\dots,\allowbreak
\Psi_{L-1}\in C(\R^{l_{L-1}},\R^{l_{L-1}})$ 
it holds that
the image of the function
\begin{equation}
  \mathcal N^{l_0,l_1,\dots,l_L}_{\Psi_1,\Psi_2,\dots,\Psi_{L-1}}
  \colon\R^{\sum_{ k = 1 }^{ L } l_k (l_{k-1} + 1)}\to C(\R^{l_0},\R^{l_L})
\end{equation}
consists precisely of the realizations of
artificial neural networks with architecture $(l_0,l_1,\dots,l_L)$ and 
activation functions $(\Psi_1,\allowbreak\Psi_2,\allowbreak\dots,\allowbreak\Psi_{L-1})$
as described in \eqref{eq:nnets} above.
In this sense the set of (realizations of) artificial neural networks
with architecture $(l_0,l_1,\dots,l_L)$ and 
activation functions $(\Psi_1,\allowbreak\Psi_2,\allowbreak\dots,\allowbreak\Psi_{L-1})$
can be para\-metrized by real vectors of length
$\sum_{ k = 1 }^{ L } l_k (l_{k-1} + 1)$,
where the parameter vector comprises
the entries of the matrices and vectors used to
describe the affine linear maps appearing in the realization function in
\eqref{eq:paramFFNN} above.
In the case where $L\in\{2,3,\dots\}$ in \cref{def:affine}
satisfies $L\geq 3$ the artificial neural network
in \eqref{eq:paramFFNN} above is referred to as a deep neural network.
One of the most prominent examples for the activation functions
$\Psi_1,\Psi_2,\dots,\Psi_{L-1}$ in \eqref{eq:paramFFNN}
is provided by the so-called multidimensional rectifier
functions (sometimes also referred to as rectified linear unit
(ReLU) functions). These rectifier functions are the subject of the next definition.
\begin{definition}[Multidimensional rectifier functions]\label{def:rect}
  Let $l\in\N$. Then we denote by
  $\rect l\colon \R^l\to\R^l$ the function
  which satisfies for all $x=(x_1,x_2,\dots,x_l)\in\R^l$
  that
  \begin{equation}
    \rect l(x)
    =
    (\max\{x_1,0\},\max\{x_2,0\},\dots,\max\{x_l,0\})
    .
  \end{equation}
\end{definition}
In the next subsection, we will employ \cref{def:nets,def:rect} above
to sketch a deep learning-based approximation scheme for linear PDEs.

\endgroup

\subsection{Towards a deep learning-based approximation method for linear PDEs}
\label{subsec:approxlinear}
We next explain how the reformulation of a linear PDE as an infinite dimensional
stochastic optimization problem from \cref{sec:reformulate} can be used
to derive a deep learning-based approximation method. For this, we now
briefly recall the setup from \cref{prop:heat_min}.
Let 
$ d \in \N $, 
$ T, \rho \in (0,\infty) $, 
$ \varrho = \sqrt{ 2 \rho T } $, 
$ a \in \R $, 
$ b \in (a,\infty) $, 
let 
$ \varphi \colon \R^d \to \R $ be a function,
let 
$
  u 
  \in C^{1,2}( [0,T] \times \R^d, \R )
$ 
be a function with at most polynomially growing partial derivatives which satisfies 
for all $ t \in [0,T] $, $ x \in \R^d $ that 
$ u(0,x) = \varphi(x) $ and 
\begin{equation}
\label{eq:differentialu2}
  \tfrac{ \partial u}{\partial t}(t,x) 
  = 
  \rho \, \Delta_x u(t,x),
\end{equation}
let 
$ (\Omega, \F, \P ) $ be a probability space, 
let
$ \mathbb{W} \colon \Omega \to \R^d $ 
be a standard normal random variable, 
let 
$ \xi \colon \Omega \to [a,b]^d $ be 
a continuous uniformly distributed random variable, 
and assume that $ \mathbb{W} $ and $ \xi $ are independent. 
\cref{prop:heat_min} then ensures 
that the solution $ u $ of the heat equation in \eqref{eq:differentialu2} 
at time $ T $ on $ [a,b]^d $ is the unique global 
minimizer of the function 
\begin{equation}
  C( [a,b]^d, \R )
  \ni v 
  \mapsto 
  \E\br[\big]{ \abs{ \varphi( \varrho \mathbb{W}  + \xi ) - v(\xi) }^2 }
  \in [0,\infty] .
\end{equation} 
Now an idea of a simple deep learning-based 
approximation method for PDEs (see~%
\cite{beck2021solving}) is 
to approximate the set $ C( [a,b]^d, \R ) $ 
of all continuous functions from $ [a,b]^d $ to $ \R $ 
through the set of all DNNs 
with a fixed sufficiently large architecture. 
More formally, 
let 
$ L \in \{ 2, 3, \dots \} $, 
$ 
  l_0, l_1, \dots, l_L
  \in \N 
$
satisfy $l_0=d$ and $l_L=1$
and consider the function\cfclear
\begin{equation}
  \begin{split}
\label{eq:DNN_search}
  \cu*{
    w \in 
    C( [a,b]^d, \R )
    \colon 
    \br*{
      \begin{gathered}
        \exists \,
        \theta \in 
        \R^{ 
          \sum_{ k = 1 }^{ L } l_k (l_{k-1} + 1)
        } 
        \colon
        \forall \, x \in [a,b]^d 
        \colon
      \\
        w(x) 
        = 
        \pr[\big]{\RealV^{l_0,l_1,\dots,l_L}_{\rect{l_1},\rect{l_2},\dots,\rect{l_{L-1}}}(\theta)}(x)
      \end{gathered}
    }
  }
  \ni v \qquad\qquad\\\qquad\qquad
  \mapsto 
  \E\br[\big]{ \abs{ \varphi( \varrho \mathbb{W}  + \xi ) - v(\xi) }^2 }
  \in [0,\infty] 
  \end{split}
\end{equation} 
\cfload.
The approach of the deep learning-based approximation algorithm in~%
\cite{beck2021solving} is then to approximately compute 
a suitable minimizer of the function 
in \eqref{eq:DNN_search} and to view the resulting approximation 
of a suitable minimizer 
of the function in \eqref{eq:DNN_search} 
as an approximation of the solution $ u $ 
of the heat equation in \eqref{eq:differentialu2} at time $ T $ on $ [a,b]^d $. 
To approximately compute a suitable minimizer of the function 
in \eqref{eq:DNN_search}, we reformulate \eqref{eq:DNN_search} by employing 
the parametrization function induced by 
the neural network description 
%in \eqref{eq:FFNN} 
in \eqref{eq:paramFFNN} 
above to obtain the function 
\begin{equation}
\label{eq:DNN_search2}
        \R^{ 
          \sum_{ k = 1 }^{ L } l_k (l_{k-1} + 1)
        } 
        \ni
        \theta 
  \mapsto 
  \E\br[\Big]{
    \abs[\big]{
      \varphi( \varrho \mathbb{W}  + \xi ) 
      - 
      \pr[\big]{\RealV^{l_0,l_1,\dots,l_L}_{\rect{l_1},\rect{l_2},\dots,\rect{l_{L-1}}}(\theta)}( \xi )
    }^2 
  }
  \in [0,\infty] .
\end{equation} 
A suitable minimizer of the function in \eqref{eq:DNN_search2} can 
then be approximately computed by means of stochastic gradient descent 
optimization algorithms.
The value of the global minimum is typically greater than $0$.

A simple implementation of the resulting algorithm in \Python\ using
the deep learning framework PyTorch is given in
\cref{lst.kolmo}.
We refer to~%
\cite{beck2021solving} for further details on such approximation schemes 
and we also refer to~%
\cite{beck2021solving} for numerical simulations for such 
approximation schemes.

\begin{center}
\begin{lstlisting}
import torch
import math

# Use the GPU if available
dev = torch.device("cuda" if torch.cuda.is_available() else "cpu")

d = 10             # the input dimension
a, b = -3.0, 3.0   # the domain will be [a,b]^d
T = 2.0            # the time horizon
ρ = 1.0            # the diffusivity

# Define the initial value
def φ(x):
    return torch.cos(x.square().sum(axis=1, keepdim=True))

# Computes an approximation of 𝔼[|φ(√(2ρT) W + ξ) - N(ξ)|²] with W a standard 
# normal random variable using the rows of x as independent realizations of the 
# random variable ξ
def loss(N, ρ, φ, T, x):
    W = torch.randn_like(x).to(dev)
    return (φ(math.sqrt(2 * ρ * T) * W + x) - N(x)).square().mean()

# Define a neural network with two hidden layers with 50 neurons each using
# ReLU activations
N = torch.nn.Sequential(
    torch.nn.Linear(d, 50), torch.nn.ReLU(),
    torch.nn.Linear(50, 50), torch.nn.ReLU(),
    torch.nn.Linear(50, 1)
).to(dev)

# Configure the training parameters and optimization algorithm
steps = 1000
batch_size = 256
optimizer = torch.optim.Adam(N.parameters())

# Train the network
for step in range(steps):
    # Generate uniformly distributed samples from [a,b]^d
    x = (torch.rand(batch_size, d) * (b-a) + a).to(dev)
    
    optimizer.zero_grad()
    # Compute the loss
    L = loss(N, ρ, φ, T, x)
    # Compute the gradients
    L.backward()
    # Apply changes to weights and biases of N
    optimizer.step()
\end{lstlisting}
  \captionof{listing}{\label{lst.kolmo}A simple implementation 
  in PyTorch of an
  algorithm based on \cref{prop:heat_min}, computing
  an approximation of the function $[-3,3]^{10}\ni x\mapsto u(2,x)\in\R$
  where $u\in C^{1,2}([0,2]\times[-3,3]^{10},\R)$ is the function 
  which satisfies for all $t\in[0,2]$, $x\in[-3,3]^{10}$ that
  $u(0,x)=\cos(\norm{x}^2)$ and
  $\tfrac{\partial u}{\partial t}(t,x) = \Delta_x u(t,x)$.
  }
  \bigskip
\end{center}

  In the above outline of a deep learning-based scheme for PDEs we have restricted 
  ourselves to the 
  numerical approximation of the simple heat equation.
  However, the procedure sketched above also extends to more general 
  linear Kolmogorov PDEs.
  In particular, the Feynman--Kac formula does not only apply to 
  linear heat PDEs as in \eqref{eq:differentialu} above,
  but can also be applied in the case of linear Kolmogorov PDEs
  with a possibly nonlinear drift coefficient function $\mu\colon\R^d\to\R^d$
  and a possibly nonlinear diffusion coefficient 
  function $\sigma\colon\R^d\to\R^{d\times d}$.
  In the PDE the function $\mu$ is then multiplied with
  the first-order partial derivatives of the solution function
  and the function $\R^d\ni x\mapsto \sigma(x)[\sigma(x)]^*\in\R^{d\times d}$
  is multiplied with the Hessian of the solution
  of the PDE; see, e.g., the PDE in (1.1) in Hairer et al.~\cite{hairer2015loss}
  and Chapter~8.1 in Øksendal \cite{oksendal2013stochastic}.
  In this case, the Feynman--Kac formula does not involve
  a simple Brownian motion as in the case of
  \cref{prop:heat_min} above,
  but instead involves a stochastic process which cannot be
  simulated exactly anymore,
  but which can be approximated by a numerical
  discretization scheme such as the 
  Euler--Maruyama or the Milstein scheme.
  We also refer to \cite{beck2021solving} 
  for further details on this more general case.

\section{Deep learning-based approximation methods for nonlinear PDEs}
\label{sec:nonlinear}

In the previous section we focused on certain deep learning-based approximation 
methods for linear heat PDEs. We will now
turn towards deep learning-based approximation methods 
for possibly nonlinear PDEs, focusing on two
particular such approximation methods and then briefly
reviewing the wider literature.
More precisely, in \cref{subsec:dgm} we will provide
a derivation of the so-called \emph{deep Galerkin method} (DGM)
proposed in Sirignano \& Spiliopoulos~\cite{Sirignano2018dgm}
(see also Dissanayake \& Phan-Thien~\cite{dissanayake1994neural}
and Lagaris et al.~\cite{lagaris1998artificial} for closely related
earlier approaches and Raissi et al.~\cite{raissi2019physics} for the
method of \emph{physics-informed neural networks} (PINNs), which
is also based on the same principle).
Then, in \cref{sec:deepsplitting}, we will give a rough outline
of the derivation of the method in \cite{beck2021deep},
which is referred to as the \emph{deep splitting} approximation method.
For both the deep Galerkin method and the deep splitting approximation
method, we will also provide simple implementations in Python using
the PyTorch framework.
We will then, in a short and necessarily very incomplete review
of the scientific literature in \cref{sec:other},
provide some brief comments on other articles
on deep learning-based approximation methods for PDEs.

\subsection{The deep Galerkin method}
\label{subsec:dgm}

The deep Galerkin method proposed in 
Sirignano \& Spiliopoulos \cite{Sirignano2018dgm}
is based on the result in \cref{thm:dgm} below, which, similar to
\cref{prop:heat_min}, casts the solution of a
certain nonlinear PDE as the solution to a
stochastic optimization problem.

For concreteness and simplicity, we state the result here only
for semilinear heat PDEs, but we note that the method can easily
be adapted to a much wider range of PDEs. We refer, e.g., to
Sirignano \& Spiliopoulos \cite{Sirignano2018dgm},
Raissi et al.~\cite{raissi2019physics}, the references
in \cref{sec:other}, and the much more comprehensive 
surveys~Cuomo et al.~\cite{cuomo2022scientific} and Karniadakis et al.~\cite{karniadakis2021physics}
for further information.

\begingroup
\providecommand{\Tv}{}
\renewcommand{\Tv}{\mathcal{T}}
\providecommand{\Xv}{}
\renewcommand{\Xv}{\mathcal{X}}
\providecommand{\Td}{}
\renewcommand{\Td}{\mathscr{t}}
\providecommand{\Xd}{}
\renewcommand{\Xd}{\mathscr{x}}

\begin{samepage}
\begin{athm}{theorem}{thm:dgm}
  Let 
    $T\in(0,\infty)$,
    $d\in\N$,
    $g\in C^2(\R^d,\R)$,
    $u\in C^{1,2}([0,T]\times\R^d,\R)$,
    $\Td\in C([0,T],(0,\infty))$,
    $\Xd\in C(\R^d,(0,\infty))$,
  assume that $g$ has at most polynomially growing partial derivatives,
  let $(\Omega,\mathcal F,\PP)$ be a probability space,
  let 
    $\Tv\colon \Omega\to[0,T]$ 
    and $\Xv\colon \Omega\to\R^d$ 
    be independent random variables,
  assume
    for all 
      $A\in\Borel([0,T])$,
      $B\in\Borel(\R^d)$
    that
    \begin{equation}
      \label{eq:densities}
      \PP(\Tv\in A)=\int_A\Td(t)\,\dd t
      \qquad\text{and}\qquad
      \PP(\Xv\in B)=\int_B \Xd(x)\,\dd x,
    \end{equation}
  let $f\colon\R\to\R$ be a Lipschitz continuous function,
  and let $\mathbb F\colon C^{1,2}([0,T]\times\R^d,\R)\to[0,\infty]$
  satisfy 
    for all
      $v=(v(t,x))_{(t,x)\in[0,T]\times\R^d}\in C^{1,2}([0,T]\times\R^d,\R)$
    that
    \begin{equation}
      \label{eq:dgmloss}
      \mathbb F(v)
      =
      \mEE[1]{
        \abs{v(0,\Xv)-g(\Xv)}^2
        +
        \abs[\big]{\tfrac{\partial v}{\partial t}(\Tv,\Xv)-\Delta_x v(\Tv,\Xv)-f(v(\Tv,\Xv))}{}^2
      }
      .
    \end{equation}
  Then the following two statements are equivalent:
  \begin{enumerate}[(i)]
    \item \llabel{it:1}
    It holds that
      $\mathbb F(u)=\inf_{v\in C^{1,2}([0,T]\times\R^d,\R)}\mathbb F(v)$.
    \item \llabel{it:2}
    It holds for all
      $t\in[0,T]$,
      $x\in\R^d$
    that
      $u(0,x)=g(x)$ and
      \begin{equation}
        \label{eq:dgm.pde}
      \tfrac{\partial u}{\partial t}(t,x) = \Delta_x u(t,x)+f(u(t,x))
      .
      \end{equation}
  \end{enumerate}
\end{athm}
\end{samepage}
\begin{aproof}
  \Nobs that 
    \cref{eq:dgmloss}
  implies that for all
    $v\in C^{1,2}([0,T]\times\R^d,\R)$
    with 
      $\forall\,x\in\R^d\colon u(0,x)=g(x)$ and
      $\forall\,t\in[0,T],\,x\in\R^d\colon \tfrac{\partial u}{\partial t}(t,x) = \Delta_x u(t,x)+f(u(t,x))$
  it holds that 
  \begin{equation}
    \llabel{eqphi}
    \mathbb F(v)=0.
  \end{equation}
    This
    and the fact that
      for all
        $v\in C^{1,2}([0,T]\times\R^d,\R)$
      it holds that
        $\mathbb F(v)\geq 0$
  establish that (\ref{thm:dgm.it:2} $\rightarrow$ \ref{thm:dgm.it:1}).
  \Nobs that
    the assumption that
      $f$ is a Lipschitz continuous function,
    the assumption that $g$ is a twice continuously differentiable function,
    and the assumption that $g$ has at most polynomially growing partial derivatives
  ensure that there exists
    $v\in C^{1,2}([0,T]\times\R^d,\R)$
  which satisfies for all
    $t\in[0,T]$,
    $x\in\R^d$
  that $v(0,x)=g(x)$ and
  \begin{equation}
    \tfrac{\partial v}{\partial t}(t,x) = \Delta_x v(t,x)+f(v(t,x))
  \end{equation}
  (cf., e.g., \cite[Corollary~3.4]{beck2021nonlinear}).
    This
    and \lref{eqphi}
  show that
  \begin{equation}
    \llabel{2}
    \inf_{v\in C^{1,2}([0,T]\times\R^d,\R)}\!\!\!\mathbb F(v) = 0.
  \end{equation}
  \Moreover \enum{
    \cref{eq:dgmloss} 
    ;
    \cref{eq:densities}
    ;
    the assumption that $\Tv$ and $\Xv$ are independent
  }[ensure] that for all
    $v\in C^{1,2}([0,T]\times\R^d,\R)$
  it holds that
  \begin{multline}
    \mathbb F(v)
    =
    \int_{[0,T]\times\R^d}
    \Bigl(\abs{v(0,x)-g(x)}^2
    \\\qquad+
    \abs[\big]{\tfrac{\partial v}{\partial t}(t,x)-\Delta_x v(t,x)-f(v(t,x))}^2\Bigr)
    \Td(t)\Xd(x)
    \,\mathrm d(t,x)
    .
  \end{multline}
    The assumption that
      $\Td$ and $\Xd$ are continuous functions
    and the fact that for all
      $t\in[0,T]$,
      $x\in\R^d$
    it holds that
      $\Td(t)\geq 0$
      and $\Xd(x)\geq 0$
  \hence imply that for all
    $v\in C^{1,2}([0,T]\times\R^d,\R)$,
    $t\in[0,T]$,
    $x\in\R^d$
    with $\mathbb F(v)=0$
  it holds that
  \begin{equation}
    \Bigl(\abs{v(0,x)-g(x)}^2
    +
    \abs[\big]{\tfrac{\partial v}{\partial t}(t,x)-\Delta_x v(t,x)-f(v(t,x))}^2\Bigr)
    \Td(t)\Xd(x)
    =
    0
    .
  \end{equation}
    This
    and the assumption that for all
      $t\in[0,T]$,
      $x\in\R^d$
    it holds that
      $\Td(t)>0$
      and $\Xd(x)>0$
  demonstrate that for all
    $v\in C^{1,2}([0,T]\times\R^d,\R)$,
    $t\in[0,T]$,
    $x\in\R^d$
    with $\mathbb F(v)=0$
  it holds that
  \begin{equation}
    \abs{v(0,x)-g(x)}^2
    +
    \abs[\big]{\tfrac{\partial v}{\partial t}(t,x)-\Delta_x v(t,x)-f(v(t,x))}^2
    =
    0
    .
  \end{equation}
  Combining 
    this 
  with
    \lref{2}
  proves that (\ref{thm:dgm.it:1} $\rightarrow$ \ref{thm:dgm.it:2}).
\end{aproof}

Similar to the approach discussed in \cref{subsec:approxlinear},
the solution to this infinite dimensional stochastic optimization
problem can be approximated by the solution to a finite dimensional
optimization problem by precomposing the functional
in \cref{eq:dgmloss} in \cref{thm:dgm} above with a function parametrizing
the set of realizations of artificial neural networks with a fixed
sufficiently large architecture and fixed activation functions.

To be more precise, under the assumptions in \cref{thm:dgm} observe that
\cref{thm:dgm} shows that the solutions of the semilinear heat PDE in
\eqref{eq:dgm.pde} are the global minimizers of the function
$\mathbb F\colon C^{1,2}([0,T]\times\R^d,\R)\to [0,\infty]$.
In a manner similar to \cref{sec:linear}, a simple deep 
learning-based approximation
method can be derived from this observation by
approximating the set $C^{1,2}([0,T]\times\R^d,\R)$
through the set of all DNNs with a fixed sufficiently large
architecture and fixed activation functions.
More formally, let $L\in\{2,3,\dots\}$,
$l_0,l_1,\dots,l_L\in\N$ satisfy
$l_0=d+1$ and $l_L=1$,
let $\Psi_1\in C^2(\R^{l_1},\R^{l_1})$,
$\Psi_2\in C^2(\R^{l_2},\R^{l_2})$,
$\dots$, $\Psi_{L-1}\in C^2(\R^{l_{L-1}},\R^{l_{L-1}})$
 and consider the function
\begin{equation}
  \label{eq:dgm.min}
  \R^{\sum_{k=1}^L l_k(l_{k-1}+1)}\ni\theta
  \mapsto
  \mathbb F\bigl(\bigl(\mathcal N^{l_0,l_1,\dots,l_L}_{\Psi_1,\Psi_2,\dots,\Psi_{L-1}}(\theta)\bigr)|_{[0,T]\times\R^d}\bigr)
  \in[0,\infty]
  .
\end{equation}
A suitable minimizer of the function in \eqref{eq:dgm.min} can then
be approximately computed using a stochastic gradient descent optimization
algorithm (where the derivatives appearing in the definition
\eqref{eq:dgmloss} of the function $\mathbb F$ are typically computed
using automatic differentiation as supplied by most popular machine learning
frameworks) and this minimizer can be considered
an approximation of a minimizer of $\mathbb F$ and thus 
an approximation of the solution
of the PDE in \eqref{eq:dgm.pde}.

\endgroup

An example implementation of the resulting algorithm is given in
\cref{lst.dgm}. For simplicity, we present the algorithm here for a
two-dimensional PDE.
In \cref{sec:simdgm} we present numerical simulations using the 
deep Galerkin method to approximate solutions of certain semilinear parabolic
PDEs in up to 200 dimensions. The source code for these
simulations can be found in \cref{sec:source_dgm}.

\begingroup
\begin{lstlisting}
import torch

# Use the GPU if available
dev = torch.device("cuda" if torch.cuda.is_available() else "cpu")

# Computes an approximation of 𝔼[|v(0,X) - φ(X)|²] using the rows of x as independent
# realizations of the random variable X
def dgm_initial_loss(v, φ, x):
    t = torch.zeros(x.shape[0],1).to(dev)
    u = v(torch.cat((t,x), axis=1))
    return (u - φ(x)).square().mean()

# Computes an approximation of 𝔼[|∂v/∂t(T,X) - Δₓv(T,X) - f(v(T,X))|²] using the 
# rows of t and x as independent realizations of the random variables T and X, 
# respectively
def dgm_dynamic_loss(v, ρ, f, t, x):
    # Split x up into components
    x1, x2 = x[:,0], x[:,1]
    # Compute v(t,x)
    u = v(torch.stack([t, x1, x2], 1))
    # Compute ∂v/∂t(t,x)
    u_t = torch.autograd.grad(u, t, torch.ones_like(u), create_graph=True)[0]
    # Compute Δₓv(t,x)
    u_x1 = torch.autograd.grad(u, x1, torch.ones_like(u), create_graph=True)[0]
    u_x1x1 = torch.autograd.grad(u_x1, x1, torch.ones_like(u_x1), 
                                 create_graph=True)[0]
    u_x2 = torch.autograd.grad(u, x2, torch.ones_like(u), create_graph=True)[0]
    u_x2x2 = torch.autograd.grad(u_x2, x2, torch.ones_like(u_x2), 
                                 create_graph=True)[0]
    Δu = u_x1x1 + u_x2x2

    return (u_t - ρ * Δu - f(u)).square().mean()

a, b = -3.0, 3.0   # the domain will be [a,b]²
T = 2.0            # the time horizon
ρ = 1.0            # the diffusivity

# The initial value of the PDE
def φ(x):
    return (1. + x.square().sum(axis=1, keepdim=True)).sqrt()

# The nonlinearity of the PDE
def f(y):
    return torch.sin(y)

# Define a fully-connected neural network with three hidden layers with
# 50 neurons in each hidden layer and using the tanh activation function
N = torch.nn.Sequential(
    torch.nn.Linear(3, 50), torch.nn.Tanh(),
    torch.nn.Linear(50, 50), torch.nn.Tanh(),
    torch.nn.Linear(50, 50), torch.nn.Tanh(),
    torch.nn.Linear(50, 1)
).to(dev)

# Configure the training parameters and optimization algorithm
steps = 250
batch_size = 256
optimizer = torch.optim.Adam(N.parameters())

for step in range(steps):
    # Generate uniformly distributed samples from [0,T]
    t = T * torch.rand(batch_size).to(dev)
    t.requires_grad_()
    # Generate uniformly distributed samples from [a,b]²
    x = (torch.rand(batch_size, 2) * (b-a) + a).to(dev)
    x.requires_grad_()
    
    optimizer.zero_grad()
    # Compute the loss
    loss = dgm_initial_loss(N, φ, x) + dgm_dynamic_loss(N, ρ, f, t, x)
    # Compute the gradients
    loss.backward()
    # Apply changes to weights and biases of N
    optimizer.step()
\end{lstlisting}
  \captionof{listing}{\label{lst.dgm}A simple implementation 
  in PyTorch of the deep Galerkin method based on \cref{thm:dgm}, computing
  an approximation of the function $u\in C^{1,2}([0,2]\times[-3,3]^{2},\R)$
  which satisfies for all $t\in[0,2]$, $x\in[-3,3]^{2}$ that
  $u(0,x)=\sqrt{1+\norm{x}^2}$ and
  $\tfrac{\partial u}{\partial t}(t,x) = \Delta_x u(t,x) + \sin(u(t,x))$.
  }
\endgroup
% \todo{Text ueberarbeiten}

\subsection{The deep splitting approximation method for semilinear PDEs}
\label{sec:deepsplitting}

We now focus on the deep splitting approximation method
introduced in \cite{beck2021deep}. The main idea of this method is
\begin{enumerate}
  \item[\mylabel{it:ds1}{($\mathcal A$)}] %
  to split up the time interval $[0,T]$, where $T\in(0,\infty)$
  is the time horizon of
  the semilinear PDE under consideration, into the subintervals 
  $[\tau_0,\tau_1]$, $[\tau_1,\tau_2]$, \ldots, $[\tau_{N-1},\tau_N]$
  where $N\in\N$, $\tau_0, \tau_1,\dots,\tau_N \in\R$ satisfy 
  $0=\tau_0<\tau_1<\dots< \tau_N = T$,
  then, 
  \item[\mylabel{it:ds2}{($\mathcal B$)}]%
  to freeze the nonlinearity of the semilinear PDE under
  consideration on each 
  of these subintervals in order to obtain a linear PDE on each
  of these subintervals,
  and, thereafter,
  \item[\mylabel{it:ds3}{($\mathcal C$)}] to consecutively apply the deep learning-based approximation method
  for linear PDEs from \cref{sec:linear} above to each subinterval.
\end{enumerate}
In this sense the deep splitting approximation method 
in \cite{beck2021deep} combines
splitting approximations (which are, in the above context, 
also referred to as exponential Euler approximations; cf.,
e.g., 
\cite{hochbruck2005explicit,cox2013pathwise,
gyongy2003splitting%
})
with deep learning-based approximations for linear PDEs (see \cref{sec:linear})
to obtain deep learning-based approximations for semilinear PDEs.
To formalize this approach, we now present in the following theorem, 
\cref{cor:splitting_with_jump_for_heat_eqn5} below,
an approximation result which illustrates how a semilinear
heat PDE can be approximated through a 
series of linear heat PDEs (according to item \ref{it:ds2} above).

\begingroup
\newcommand{\f}{f}
\renewcommand{\u}{u}
\newcommand{\U}[3]{\mc U^{#1,#2}_{#3}}

\begin{athm}{theorem}{cor:splitting_with_jump_for_heat_eqn5}
  Let 
  $T\in(0,\infty)$,
  $p\in[1,\infty)$,
  $f\in C^2(\R,\R)$,
  let $\u_d\in C^{1,2}([0,T]\times\R^d,\R)$, $d\in\N$, satisfy for all
  $d\in\N$,
	$t\in[0,T]$,
	$x\in\R^d$
  that
  \begin{equation}
    \label{eq:semilinearheat}
	(\tfrac{\partial}{\partial t}\u_d)(t,x)
	=
	(\Delta_x \u_d)(t,x)+\f(\u_d(t,x))
	,
	\end{equation}
  and assume for all
  $d\in\N$,
	$i,j\in\{1,2,\ldots,d\}$ 
	that 
  $\sup_{t\in[0,T]}\sup_{x=(x_1,x_2,\dots,x_d)\in\R^d}
  \br[\big]{(1\allowbreak+\sum_{k=1}^d\abs{x_k})^{-p}\allowbreak
  (
  \abs{(\tfrac{\partial^2}{\partial x_i\partial x_j}u_d)(t,x)}+\abs{(\tfrac{\partial}{\partial t}u_d)(t,x)}+\abs{f''(x_1)})+\abs{f'(x_1)}
  }<\infty$.
  Then 
  \begin{enumerate}[(i)]
    \item \label{item:nonlinear1}
    there exist unique at most polynomially growing
      $\U dNn \in C^{1,2}([\tfrac{(n-1)T}{N},\tfrac{nT}{N}]\times\R^d,\R)$, $d,N\in\N$, $n\in\{0,1,\ldots,N\}$,
    which satisfy
      for all
        $d,N\in\N$,
        $n\in\{0,1,\dots,N-1\}$,
        $t\in[\frac{nT}N,\frac{(n+1)T}N]$,
        $s\in[\frac{-T}N,0]$,
        $x\in\R^d$
      that
        $\U dN{n+1}(\tfrac{nT}{N},x)=\U dNn(\tfrac{nT}{N},x)+\tfrac{T}N\f(\U dNn(\tfrac{nT}{N},x))$,
        $\U dN0(s,x)=\u_d(0,x)$,
        and
        \begin{flexmath}d
          \label{eq:linearheat}
          \pr[\big]{\tfrac{\partial}{\partial t}\U dN{n+1}}(t,x)
          =
          \pr[\big]{\Delta_x\U dN{n+1}}(t,x)
        \end{flexmath}
    and
    \item \label{item:nonlinear2} 
    there exists
      $c\in\R$
    such that for all
      $d,N\in\N$,
      $x=(x_1,x_2,\dots,x_d)\in\R^d$
    it holds that
    \begin{flexmath}d[eq:claim]
      \abs{\U dNN(T,x)-\u_d(T,x)}
      \\
      &\leq 
      cd^{p+1}N^{-\nicefrac12}\pr[\big]{1+\textstyle\sum_{i=1}^d\abs{x_i}}^p
      .
    \end{flexmath}
  \end{enumerate}
\end{athm}%
We now add some comments on the mathematical objects
appearing in \cref{cor:splitting_with_jump_for_heat_eqn5}
above. 
The semilinear heat PDEs in \eqref{eq:semilinearheat} describe
the PDEs whose solutions we intend to approximate in
\cref{cor:splitting_with_jump_for_heat_eqn5}.
The real number $T\in(0,\infty)$ 
describes the 
time horizon of the semilinear heat PDEs
in \eqref{eq:semilinearheat}.
The function $f\colon \R\to\R$ 
in \cref{cor:splitting_with_jump_for_heat_eqn5}
specifies the nonlinearity
of the semilinear heat PDEs in \eqref{eq:semilinearheat}.
The functions $u_d\colon[0,T]\times\R^d\to\R$, 
$d\in\N$, 
in \cref{cor:splitting_with_jump_for_heat_eqn5}
describe the exact solutions of the
semilinear heat PDEs in \eqref{eq:semilinearheat}.

In \cref{cor:splitting_with_jump_for_heat_eqn5} 
we assume that the solutions 
$u_d\colon[0,T]\times\R^d\to\R$, 
$d\in\N$, of the semilinear heat PDEs
in \eqref{eq:semilinearheat} grow at most polynomially.
The real number $p\in[1,\infty)$
is used to formulate this polynomial growth
assumption in \cref{cor:splitting_with_jump_for_heat_eqn5}.
More formally, in \cref{cor:splitting_with_jump_for_heat_eqn5}
it is assumed that for all $d\in\N$ there
exists $\mathfrak c\in\R$ such that for all
$t\in[0,T]$,
$x=(x_1,x_2,\dots,x_d)\in\R^d$ 
it holds that 
$\abs{u_d(t,x)}\leq \mathfrak c(1+\sum_{k=1}^d\abs{x_i})^p$.
In \cref{cor:splitting_with_jump_for_heat_eqn5} we also assume that
$f\colon\R\to\R$ is Lipschitz continuous and that
the first derivatives of $u_d\colon[0,T]\times\R^d\to\R$, $d\in\N$, the second derivatives
of $u_d\colon[0,T]\times\R^d\to\R$, $d\in\N$, and the second derivative of
$f\colon\R\to\R$ grow at most polynomially.
More specifically, the assumption in \cref{cor:splitting_with_jump_for_heat_eqn5} 
that for all
$d\in\N$, $i,j\in\{1,2,\dots,d\}$ 
it holds that 
\begin{multline}
  \sup_{t\in[0,T]}\sup_{x=(x_1,x_2,\dots,x_d)\in\R^d}\br[\Big]{\pr[\big]{1+{\textstyle\sum}_{k=1}^d\abs{x_k}}^{-p}\pr[\Big]{\abs{u_d(t,x)}+\abs[\big]{\pr[\big]{\tfrac{\partial}{\partial x_i}u_d}(t,x)}
  \\+\abs[\big]{\pr[\big]{\tfrac{\partial^2}{\partial x_i\partial x_j}u_d}(t,x)}+\abs[\big]{\pr[\big]{\tfrac{\partial}{\partial t}u_d}(t,x)}+\abs{f''(x_1)}}+\abs{f'(x_1)}}<\infty 
\end{multline}
ensures that for all
$d\in\N$ there exists $\kappa\in\R$ such that for all
$i,j\in\{1,2,\dots,d\}$,
$t\in[0,T]$,
$x=(x_1,x_2,\dots,x_d)\in\R^d$,
$v,w\in\R$
it holds that
$\abs{f(v)-f(w)}\leq \kappa\abs{v-w}$,
$\abs{(\tfrac{\partial}{\partial x_i}u_d)(t,x)}\leq \kappa(1+\sum_{k=1}^d\abs{x_i})^p$,
$\abs{(\tfrac{\partial^2}{\partial x_i\partial x_j}u_d)(t,x)}\leq \kappa(1+\sum_{k=1}^d\abs{x_i})^p$,
and $\abs{f''(v)}\leq\kappa(1+\abs{v})^p$.

\Cref{cor:splitting_with_jump_for_heat_eqn5}
shows that on the subintervals 
$\br[\big]{0,\tfrac TN}$, $\br[\big]{\tfrac TN,\tfrac{2T}N}$, \ldots, $\br[\big]{\tfrac{(N-1)T}N,T}$
for $N\in\N$ 
there exist unique at most polynomially
growing classical solutions $\mathcal U_N^{d,n}\colon\br[\big]{\tfrac{(n-1)T}N,\tfrac{nT}N}\times\R^d\to\R$, $d,N\in\N$,
$n\in\{1,2,\dots,N\}$, of the
linear heat PDEs in \eqref{eq:linearheat} whose initial conditions
are specified by means of the initial conditions
$\R^d\ni x\mapsto u_d(0,x)\in\R$, $d\in\N$, 
of the solutions $u_d\colon [0,T]\times\R^d\to\R$,
$d\in\N$, of the semilinear heat PDEs
in \eqref{eq:semilinearheat}.
\Cref{item:nonlinear2} in \cref{cor:splitting_with_jump_for_heat_eqn5}
ensures that these unique at most polynomially growing classical solutions
$\U dNn\colon\br[\big]{\tfrac{(n-1)T}N,\tfrac{nT}N}\times\R^d\to\R$, $d,N\in\N$,
$n\in\{1,2,\dots,N\}$, of the linear heat PDEs in
\eqref{eq:linearheat} can then be used to approximate the solutions
$u_d\colon[0,T]\times\R^d\to\R$, $d\in\N$, of the semilinear heat PDEs in \eqref{eq:semilinearheat}.
In particular, observe that \cref{item:nonlinear2} in \cref{cor:splitting_with_jump_for_heat_eqn5}
proves that for every $d\in\N$, $x\in\R^d$
it holds that the approximation error 
$\abs{\U dNN(T,x)-\u_d(T,x)}$ converges to $0$ as the number of subintervals
$N$ goes to $\infty$.
\endgroup%

\Cref{cor:splitting_with_jump_for_heat_eqn5} can be proved through 
applications of the Feynman--Kac formula for linear heat equations 
(cf., e.g., \cite[Corollary~4.17]{hairer2015loss}), the Feynman--Kac formula for semilinear
heat equations (cf., e.g., \cite[Theorem~1.1]{beck2021nonlinear}),
and the Gronwall inequality.
The proof of \cref{cor:splitting_with_jump_for_heat_eqn5} is 
essentially standard, and in order to keep this overview article at a reasonable length
we omit here the detailed proof of \cref{cor:splitting_with_jump_for_heat_eqn5}
(cf., e.g.,~Germain et al.~\cite[Section~3.2]{GermainEtal2020}, Cox \& van Neerven~\cite[Section~3]{cox2013pathwise},
and Jentzen~\cite[Theorem~2]{jentzen2009pathwise2}).

We present in 
\cref{lst.ds} below a simple implementation of the algorithm sketched in 
\cref{it:ds1,it:ds2,it:ds3} above based on \cref{cor:splitting_with_jump_for_heat_eqn5} above.
In \cref{sec:simds} we present numerical simulations approximating certain semilinear parabolic PDEs
in up to 1000 dimensions. The source code for these simulations can be found in \cref{sec:source_ds}.
We refer to~%
\cite{beck2021solving} for further details on the deep splitting approximation scheme
and we also refer to~%
\cite{beck2021solving} for further numerical simulations for this scheme.

\begingroup
\begin{lstlisting}
import torch
import math

# Use the GPU if available
dev = torch.device("cuda" if torch.cuda.is_available() else "cpu")

# Computes an approximation of 𝔼[|φ(√(2ρT) W + ξ) - N(ξ)|²] with W a standard 
# normal random variable using the rows of x as independent realizations of the 
# random variable ξ
def loss(N, ρ, φ, T, x):
    W = torch.randn_like(x).to(dev)
    return (φ(math.sqrt(2 * ρ * T) * W + x) - N(x)).square().mean()

d = 10             # the input dimension
a, b = -3.0, 3.0   # the domain will be [a,b]^d
T = 2.0            # the time horizon
ρ = 1.0            # the diffusivity

# The initial value of the PDE
def φ(x):
    return torch.cos(x.square().sum(axis=1, keepdim=True))

# The nonlinearity of the PDE
def f(y):
    return torch.sin(y)

n = 30             # the number of time subintervals

# Configure the training parameters
steps = 1000
batch_size = 256

# Define n neural networks with two hidden layers with 50 neurons each
# using ReLU activations
N = [
    torch.nn.Sequential(
        torch.nn.Linear(d, 50), torch.nn.ReLU(),
        torch.nn.Linear(50, 50), torch.nn.ReLU(),
        torch.nn.Linear(50, 1)
    ).to(dev) for _ in range(n)
]

for i in range(n):
    # In this iteration we train the network N[i]
    optimizer = torch.optim.Adam(N[i].parameters())

    # Define the initial value for the linear PDE problem
    # on the subinterval [iT/n,(i+1)T/n]
    def ψ(x):
        y = φ(x) if i==0 else N[i-1](x)
        return y + T/n * f(y)

    # Train the network
    for step in range(steps):
        # Generate uniformly distributed samples from [a,b]^d
        x = (torch.rand(batch_size, d) * (b-a) + a).to(dev)
        
        optimizer.zero_grad()
        # Compute the loss
        L = loss(N[i], ρ, ψ, T/n, x)
        # Compute the gradients
        L.backward()
        # Apply changes to weights and biases of N[i]
        optimizer.step()
\end{lstlisting}
  \captionof{listing}{\label{lst.ds}A simple implementation 
  in PyTorch of the deep splitting method associated to \cref{cor:splitting_with_jump_for_heat_eqn5}, 
  computing for each $i\in\{1,2,\dots,30\}$ an approximation 
  of the function $[-3,3]^{10}\ni x\mapsto u(\frac{i}{15},x)\in\R$
  where $u\in C^{1,2}([0,2]\times[-3,3]^{10},\R)$ is the function which
  satisfies for all $t\in[0,2]$, $x\in[-3,3]^{10}$ that
  $u(0,x)=\cos(\norm{x}^2)$ and
  $\tfrac{\partial u}{\partial t}(t,x) = \Delta_x u(t,x) + \sin(u(t,x))$.
  }
\endgroup

\subsection{Other deep learning-based approximation methods for PDEs}
\label{sec:other}

The methods sketched in \cref{subsec:approxlinear,subsec:dgm,sec:deepsplitting}
are just three of a large number of machine learning-based approximation schemes for PDEs 
that have been proposed in the scientific literature. In the following, we provide a
short overview of the literature in some strands of research in this area. We make no
claims whatsoever to comprehensiveness; in fact, due to the sheer size and rapid development
of the field, only a very incomplete view can be attempted here.

\subsubsection*{Approximation methods for linear PDEs}
We refer, e.g., to 
\cite{beck2021solving} for approximation
methods for linear Kolmogorov PDEs based on discretizations
of stochastic differential equations.
We refer, e.g., to Sabate Vidales et al.~\cite{VidalesEtal2019}
for approximation methods for linear Kolmogorov PDEs
based on
discretizations of stochastic differential equations 
in conjunction with suitable control variates.
We refer, e.g., to Khoo \& Ying \cite{khoo2019switchnet} for approximation
methods for scattering problems associated with linear PDEs
of the Helmholtz type
based on representing the forward and inverse map by DNNs.

\subsubsection*{Approximation methods for PDEs based on
formulations in terms of backward stochastic differential equations (BSDEs)}

The so-called \emph{deep BSDE} method is
based on exploiting the equivalence between 
certain PDEs and BSDEs
and approximating the gradient of the solution at discrete time points
by DNNs (cf., e.g., Pardoux \& Peng~\cite{pardoux1990adapted} for the connection between
decoupled BSDEs and semilinear PDEs and
cf., e.g., \cite{pardoux1992backward,peng1991probabilistic} for the connection between
coupled BSDEs and quasilinear PDEs). This method was first introduced 
in E et al.~\cite{EHanJentzen2017} and
Han et al.~\cite{HanJentzenE2018} for semilinear parabolic PDEs.

Related approaches for semilinear parabolic PDEs can be found, e.g., in
Hur\'e et al.~\cite{HurePhamWarin2019} and Chan-Wai-Nam et al.~\cite{ChanMikaelWarin2019},
with both these references also containing extensive numerical simulations
comparing different approaches. We refer, e.g., to Henry~\cite{Henry2017deep} and 
Pereira et al.~\cite{WangEtal2019} for
approximation methods for BSDEs building on the deep BSDE approach applied to problems from
mathematical finance and robotics, respectively.
We refer, e.g., to Han et al.~\cite{han2020solving} for 
methods related to the deep BSDE approach suitable for
the approximate solution of eigenvalue problems for 
semilinear second order differential operators;
we refer, e.g., to Fujii et al.~\cite{fujii2019asymptotic} for refinements of 
the deep BSDE method %
suitable for the approximative pricing of American options;
we refer, e.g., to G\"uler et al.~\cite{GulerEtal2019} 
for improvements and numerical results
for such approximation methods;
and we refer, e.g., to Castro~\cite{castro2021deep} 
for extensions of the method introduced in 
\cite{HurePhamWarin2019}
aimed towards solving non-local nonlinear PDEs.
We also refer, e.g., to N\"usken \& Richter~\cite{nusken2020solving} for
approximation methods for Hamilton--Jacobi--Bellman
PDEs and control problems based on discretizations
of decoupled BSDEs.

The connection between quasilinear parabolic PDEs and coupled BSDEs
is exploited, e.g., in Raissi~\cite{raissi2018forward}
and Han \& Long~\cite{HanLong2020} to obtain deep learning-based
approximation schemes for quasilinear parabolic PDEs.
We refer, e.g., to Ji et al.~\cite{JiEtal2020},
Andersson et al.~\cite{andersson2022convergence},
and Jiang \& Li \cite{jiang2021convergence}
for further approximation methods for coupled BSDEs.
We refer to Gonon et al.~\cite{GononEtal2020} for a study of 
the approximation methods in \cite{HanJentzenE2018,HanLong2020}
when applied to problems of finding risk-sharing equilibria
for asset pricing.

In Kremsner et al.~\cite{kremsner2020deep} the deep BSDE
method has been adapted to the approximation of semilinear elliptic PDEs
based on time discretizations of BSDEs with random terminal times.
Approximation methods based on BSDE representations have also been
adapted to fully nonlinear PDEs. We refer, e.g., to Pham et al.~\cite{pham2020neural}
and Germain et al.~\cite{GermainEtal2020} for extensions of the approach in
\cite{HurePhamWarin2019} to fully nonlinear PDEs and we refer, e.g., to
Beck et al.~\cite{beck2019machine}
and Pereira et al.~\cite{PereiraEtal2020}
for approximation methods for
fully nonlinear PDEs based on extensions of the deep BSDE method
employing time discretizations of second order BSDEs.
We also refer to Jacquier \& Oumgari~\cite{jacquier2019deep} for
approximation methods for certain path-dependent PDEs
based on discretizations of BSDEs.

\subsubsection*{Approximation methods based on least square 
minimization of the residual}
Some of the first attempts at approximatively solving PDEs using
artificial neural networks were based on an idea similar to the one
presented in \cref{subsec:dgm} above, i.e., least squares minimization
of the residuals of the PDEs under consideration; see, e.g.,
Dissanayake \& Phan-Thien~\cite{dissanayake1994neural}
and Lagaris et al.~\cite{lagaris1998artificial} for
early articles using this approach. These works use a fixed deterministic
set of collocation points and are hence only suitable for
low-dimensional PDEs. The \emph{physics-informed neural networks} (PINNs)
introduced by Raissi et al.~\cite{raissi2017physics1,raissi2017physics2} 
(published as \cite{raissi2019physics}) employ a
similar approach but using a fixed set of randomly chosen collocation
points throughout training (and integrating observed
data into the learning process), see also 
Berg \& Nystr\"om~\cite{berg2018unified}
for a similar early example of this approach
suitable for complex domains.
The DGM variant of the least-squares minimization method presented in
\cref{subsec:dgm} above, based on variable stochastic collocation
points was first introduced in Sirignano \& Spiliopoulos~\cite{Sirignano2018dgm}
(see also Carleo \& Troyer~\cite{carleo2017solving} for a similar method
using restricted Boltzmann machines in the approximate solution
of Schr\"odinger PDEs with spins).

Many interesting extensions and variations of the 
PINN/DGM approach have been proposed in the scientific literature.
We refer, e.g., to Anitescu et al.~\cite{AnitescuEtal2019} 
for approximation methods
for second order boundary value PDE problems
based on adaptive deterministic collocation points,
Wang \& Perdikaris~\cite{WangPerdikaris2020}
treats approximation methods for free boundary and Stefan PDE problems.
In Yang et al.~\cite{yang2020physics} a PINN approach and generative adversarial
networks are combined to obtain approximation methods for 
stochastic differential equations.
We also refer, e.g., to Zhang et al.~\cite{zhang2020learning} for approximation methods 
employing PINNs for low-dimensional nonlinear stochastic PDEs based on
generalized Karhunen--Lo\`eve type
expansions.

The PINN/DGM approach has also been extended to integral equations
and non-local PDEs. Some of these methods use classical techniques
for the non-local part (see, e.g., Pang et al.~\cite{pang2019fpinns},
Pang et al.~\cite{PANG2020109760},
Lu et al.~\cite{LuEtal2020},
and Mishra \& Molinaro \cite{MISHRA2021107705}) 
and are hence
susceptible to the curse of dimensionality, but approaches
that are suitable for approximating solutions of high-dimensional PDEs
have also been recently investigated (see, e.g., 
Al-Aradi et al.~\cite{alaradi2019extensions},
Deveney et al.~\cite{DeveneyEtal2020},
Guo et al.~\cite{guo2022monte},
and Yuan et al.\cite{yuan2022apinn}).
We also refer to Haghighat et al.~\cite{haghighat2021nonlocal} for related
attempts at incorporating non-local interactions into
a PINN/DGM setting.

There is also a large amount of research focusing on the choice of
architectures and methods for improving the training in schemes involving
the minimization of the PDE residual;
for examples of such works, we refer, e.g., to Meng et al.~\cite{meng2020ppinn}
and Jagtap et al.~\cite{Jagtap20202002,JAGTAP2020113028}
for methods splitting up the time and/or space domain 
in order to improve and parallelize training (see also Hu et al.~\cite{hu2022when} for
numerical results and analyses of the generalization error for one such method);
we refer, e.g., to Yu et al.~\cite{yu2022gradient},
McClenny \& Braga-Neto~\cite{mcclenny2020self},
Wight \& Zhao~\cite{wight2020solving}, and Xu et al.~\cite{xu2022adaptive}
for adaptive sampling strategies;
we refer, e.g., to Jagtap et al.~\cite{jagtap2020locally,JAGTAP2020109136}
for training strategies involving activation functions with learnable parameters;
we refer, e.g., to Wang et al.~\cite{WANG2021113938},
Wang et al.~\cite{wang2021understanding},
and Li et al.~\cite{li2022}
for proposals of particular architectures suitable for PINNs;
we refer, e.g., to Krishnapriyan et al.~\cite{krishnapriyan2021characterizing}
and Li et al.~\cite{li2022} for strategies involving weighting of the 
various terms in the loss function;
and we refer, e.g., to Yu et al.~\cite{yu2022gradient}
for methods involving including the gradient of the residual in the loss function.

For theoretical results giving bounds on the generalization error in
the training of DNNs in schemes involving least-squares minimzation of the
residual of PDEs, we refer, e.g.,
to Shin et al.~\cite{shin2020convergence},
Mishra \& Molinaro~\cite{mishra2021estimates},
and Hu et al.~\cite{hu2022when}.

Numerous articles in the literature focus on analyzing the performance of
the \mbox{PINN/DGM} approach for particular applied problems, see, e.g.,
Kissas et al.~\cite{KISSAS2020112623},
Cai et al.~\cite{Cai2021},
Mao et al.~\cite{MAO2020112789},
Cai et al.~\cite{cai2021physics},
and Jin et al.~\cite{JIN2021109951} (though note that these articles
treat low-dimensional problems).
For further numerical studies of approximation
methods for PDEs based on the PINN/DGM approach,
we refer, e.g., to
Dockhorn~\cite{dockhorn2019discussion},
Michoski et al.~\cite{MichoskiEtal2020},
Magill et al.~\cite{magill2018neural},
and Gorikhovskii et al.~\cite{Gorikhovskii_2022}.

Though we hope that the above references provide a useful introduction to
various directions of research within this area, we note that the PINN/DGM approach has spawned an
enormous amount of research activity, which we could not possibly survey
in any comprehensive manner here. Cuomo et al.~\cite{cuomo2022scientific}
provides a recent overview of this field (see also 
Karniadakis et al.~\cite{karniadakis2021physics}).

Finally, we refer to Blechschmidt \& Ernst~\cite{blechschmidt2021three} for
a presentation and comparison of three deep learning-based approaches to
approximatively solving PDEs, including
one based on least-squares minimization of the residual as well as a BSDE-based approch,
notably accompanied by annotated source 
code in the form of Jupyter notebooks.

\subsubsection*{Other approximation methods}

For extensions of the deep splitting method as presented in
\cref{sec:deepsplitting} above, we refer to \cite{beck2020deep}
for approximation
methods for possibly high-dimensional semilinear stochastic PDEs, 
and we refer to Boussange et al.~\cite{boussange2022deep}
and Frey \& K\"ock \cite{frey2022deep}
for approximation methods for non-local
PDEs with various boundary conditions.

In E \& Yu~\cite{E2018deep}
the so-called \emph{deep Ritz} method, an approximation method for elliptic PDEs
based on suitable variational formulations for the
PDEs under consideration,
is introduced, see also the related methods in, e.g., 
Khoo et al.~\cite{khoo2019solving}
and Wang \& Zhang~\cite{WangZhang2020}
and theoretical results in M\"uller \& Zeinhofer~\cite{muller2020deep}
and Karumuri et al.~\cite{KARUMURI2020109120}. 
Extensions of such methods suitable for
boundary value problems on complex geometries
can be found, e.g., in
Li et al.~\cite{li2020d3m}
and Sheng \& Yang~\cite{SHENG2021110085}.
In Samaniego et al.~\cite{samaniego2020energy} a similar method
is applied to problems from computational mechanics.
Comparison of methods based on a variational formulation with ones
based on least-square minimzation of the residual can be found, e.g., in 
Nabian \& Mohammad \cite{NABIAN201914}
and Karumuri et al.~\cite{KARUMURI2020109120}.

Approximation methods based on a Petrov--Galerkin formulation
(mostly suitable for low-dimensional PDEs) can be found, e.g.,
in Kharazmi et al.~\cite{kharazmi2019variational}.
Refinements of this method can be found, e.g., in
Kharazmi et al.~\cite{KHARAZMI2021113547}, for further numerical results
and error analyses we refer, e.g., to Berrone et al.~\cite{Berrone2022,Berrone2022b}.

We refer, e.g., to Han et al.~\cite{HanEtal2020} for approximation
methods for quasilinear elliptic PDEs
based on suitable Feynman--Kac type representations for quasilinear
elliptic PDEs.
We refer, e.g., to Nakamura-Zimmerer et al.~\cite{NakamuraZimmerer2020} for approximation
methods for Hamilton--Jacobi--Bellman PDEs based on
Pontryagin's Minimum Principle and adaptive 
sampling strategies.
We refer, e.g., to Raissi et al.~\cite{raissi2018numerical} for
approximation methods for nonlinear PDEs based on
Gaussian process regression.

We refer, e.g., to Zang et al.~\cite{Zang2020weak} for approximation
methods for nonlinear PDEs based on approximating
weak solutions of PDEs using adversarial neural networks.

\subsubsection*{Particular problems}

We refer, e.g., to \cite{BeckerCheridito2019,BeckerCheriditoJentzen2019,
HurePhamWarin2019,ChenWan2020,BayerEtal2020,
goudenege2019variance}
for approximation methods for American option problems
associated with free boundary PDE problems
(cf., e.g., \cite{vanmoerbeke1976optimal},  \cite[Chapter~VII]{peskir2006optimal},
and \cite[Section~12.3]{oksendal2013stochastic}
for the connection between American options
and free boundary PDE problems).
We refer, e.g., to \cite{han2019solving,Hermann2019deep,
luo2019backflow,pfau2019ab,cai2018approximating}
for approximation methods for many-electron Schrödinger
PDEs.
We refer, e.g., to \cite{uchiyama1993solving,chen2020physics,fan2019solving,
fan2020solving,khoo2019switchnet,
raissi2019physics,zhang2019quantifying,
raissi2018deep,LongLuMaDong2018,LongLuDong2019,LuEtal2020,
WangPerdikaris2020,DarbonLangloisMeng2020}
for approximation methods for solving inverse problems
associated to PDEs.
We refer, e.g., to Han \& Hu~\cite{Han2019deep} 
and Han et al.~\cite{han2020convergence}
for approximation methods for Markovian Nash equilibria
of stochastic games with a finite number of agents
based on approximations for Hamilton--Jacobi--Bellman PDEs
and fictitious play.
We refer, e.g., to \cite{Carmona2019convergence1,
Carmona2019convergence2,Lin2020apac,Ruthotto2020machine}
for approximation methods for Markovian Nash equilibria
of stochastic games with an infinite number of agents
based on mean-field game theory and
approximations for Hamilton--Jacobi--Bellman and
Fokker--Planck PDEs.
We refer, e.g., to Lye et al.~\cite{lye2020deep} for approximation
methods for observables of 
parameterized convection-diffusion PDEs.
We refer, e.g., to 
\cite{khoo2019solving,khoo2020solving,ZhuEtal2019,LongoEtal2020,LyeEtal2020,DeveneyEtal2020,brandstetter2022clifford}
for approximation methods for parametric PDE problems
based on surrogate modelling.

For further approximation methods 
for parametric PDEs and corresponding numerical results
we refer, e.g., to
Geist et al.~\cite{GeistEtal2020},
Li et al.~\cite{li2020fourier} and Becker et al.~\cite{becker2022learning}.

\subsubsection*{Software libraries}
For software libraries implementing deep learning-based approximation
methods for PDEs, we refer, e.g., to the Python-based libraries DeepXDE
(\url{deepxde.readthedocs.io}; cf.~Lu et al.~\cite{LuEtal2020})
and
neurodiffeq (\url{pypi.org/project/neurodiffeq}; cf.~Chen et al.~\cite{chen2020neurodiffeq})
both implementing the \mbox{PINN/DGM} method
and the Julia-based libraries 
NeuralPDE.jl (\url{neuralpde.sciml.ai}; cf.~Zubov et al.~\cite{zubov2021neuralpde})
implementing the PINN/DGM method as well as the deep BSDE method
and HighDimPDE.jl (\url{highdimpde.sciml.ai}; cf.~Boussange et al.~\cite{boussange2022deep}),
implementing the deep splitting and full history recursive
multilevel Picard (see \cref{sec:othermethods} below) methods.

\section[Non-machine learning-based approximation methods for PDEs]{Approximation methods for high-dimensional PDEs not based on machine learning}
\label{sec:othermethods}
Beside machine learning-based approximation methods 
for PDEs,
there are also several other attempts in the scientific literature to approximately
solve high-di\-men\-sional nonlinear PDEs.
In particular, we refer, e.g., to Darbon \& Osher~\cite{DarbonOsher2016}
for approximation methods for certain high-dimensional
first-order Hamilton--Jacobi--Bellman PDEs.
We refer, e.g., to~\cite{BallyPages2003a,%
BallyPages2003b,%
DelarueMenozzi2006,%
DelarueMenozzi2008,%
zhao2014new,
fu2016multistep,
teng2020multi,
zhao2010stable
} for approximation methods for BSDEs based on spatial grids.
We refer, e.g., to~\cite{%
BouchardTouzi2004,%
CrisanManolarakisTouzi2010,%
GobetTurkedjiev2016Malliavin, %
hu2011malliavin
} for approximation methods for BSDEs based on Malliavin calculus.
We refer, e.g., to~\cite{%
GobetLemor2008,%
GobetLemorWarin2005,%
LemorGobetWarin2006,%
bender2012least,%
gobet2016stratified,%
GobetTurkedjiev2016regression, %
RuijterOosterlee2015,%
RuijterOosterlee2016%
} for approximation methods for BSDEs based on suitable projections on function spaces.
We refer, e.g., to~\cite{
chassagneux2020cubature,%
CrisanManolarakis2012,%
CrisanManolarakis2014,%
de2015cubature
} for approximation methods for BSDEs based on cubature on Wiener space.
We refer, e.g., to~\cite{%
chang2016branching,
le2017particle,
le2018monte,
LeCavilOudjaneRusso2019%
} for approximation methods for PDEs based on density estimations and particle systems.
We refer, e.g., to~\cite{%
BenderDenk2007,%
bender2008time,%
GobetLabart2010,%
labart2013parallel%
} for approximation methods based on Picard iterations and suitable projections on function spaces.
We refer, e.g., to~\cite{%
BriandLabart2014,%
GeissLabart2016
} for approximation methods for BSDEs based on Wiener chaos expansions.
We refer, e.g., to~\cite{%
briand2001donsker,
ma2002numerical,
martinez2011numerical,
cheridito2013bs,
geiss2020random,
geiss2020mean
} for approximation methods for BSDEs based on random walks.
We refer, e.g., to~\cite{%
chassagneux2015numerical,
imkeller2010results,
richou2011numerical,
richou2012markovian
} for approximation methods for quadratic BSDEs.
We refer, e.g., to Abbas-Turki et al.~\cite{%
abbasturki2020conditional,
abbas2020conditional
} for approximation methods for BSDEs based on nested Monte Carlo approximations.
We refer, e.g., to~\cite{skorokhod1964branching,watanabe1965branching,Henry-Labordere2012,
Henry-Labordere2014,chang2016branching,warin2017variations,henry2019branching}
for approximation methods for semilinear parabolic PDEs based on branching diffusion
approximations.
We refer, e.g., to~\cite{richter2021solving,dolgov2021tensor,eigel2017adaptive,horowitz2014linear}
for approximation methods for BSDEs as well as
parabolic and Hamilton--Jacobi--Bellman PDEs based on
tensor trains.
We refer, e.g., to Warin~\cite{warin2018nesting,warin2018monte%
}
for approximation methods for semilinear parabolic PDEs based on 
standard Monte Carlo approximations for nested conditional expectations.

\subsubsection*{Full history multilevel Picard approximation methods}

It should be noted that despite the large amount
of research and the great potential that
deep-learning based approximation methods for high-dimensional
PDEs have shown in numerical simulations,
as of today, it has not been rigorously proved that they indeed
overcome the curse of dimensionality in the approximation
of high-dimensional PDEs. To our knowledge, the only
approximation methods for which it has been established
in the scientific literature
that they overcome the curse of dimensionality in the 
numerical approximation
of semilinear PDEs with general time horizons
make up the class of 
full history recursive multilevel Picard approximation
methods (in the following we will abbreviate
\emph{full history recursive multilevel Picard} by MLP).
MLP approximation methods were introduced in
2016 in E et al.~\cite{E2016multilevel,E2019multilevel} and
in 2018 in Hutzenthaler et al.~\cite{hutzenthaler2020overcoming}
and have been extended and further studied analytically and numerically in
\cite{hutzenthaler2020multilevel,%
hutzenthaler2019arxiv1903,beck2020overcoming,%
giles2019generalised,hutzenthaler2019overcoming,beck2020arxiv2003,%
becker2020numerical,hutzenthaler2020lipschitz,beck2020nested}.
Roughly speaking, MLP approximation methods are based on the idea, first,
\begin{enumerate*}[(I)]
  \item
  to reformulate the PDE problem under
consideration as a suitable stochastic fixed point
equation, then,
  \item
  to approximate the fixed point of the resulting
stochastic fixed point equation through fixed point
iterates, which in the context of temporal integral
equations are referred to as Picard iterations,
and, finally,
\item to approximate the expectations and the integrals
appearing in the fixed point iterates
through suitable multilevel Monte Carlo approximations.
\end{enumerate*}
The resulting approximations are full history recursive in the
sense that the calculation of the $n$th MLP iterate employs realizations
of the $(n-1)$th, $(n-2)$th, \ldots, $1$st MLP iterate.
The full history recursive nature of MLP approximation methods
is also one of the major differences of MLP approximation
methods when compared to standard multilevel Monte Carlo
approximations (for references on standard multilevel Monte Carlo
approximations see, e.g., 
\cite{heinrich1998monte,heinrich1999monte,heinrich2001multilevel,%
giles2008multilevel,giles2015multilevel}).

We refer to~\cite{E2016multilevel,hutzenthaler2020overcoming,
giles2019generalised}
for MLP approximation methods
for semilinear heat PDEs.
We refer to~\cite{%
hutzenthaler2020multilevel,
hutzenthaler2019overcoming}
for MLP approximation methods
for semilinear heat PDEs with gradient-dependent nonlinearities.
We refer to~\cite{%
hutzenthaler2019arxiv1903,
hutzenthaler2020lipschitz}
for MLP approximation methods
for semilinear parabolic PDEs.
We refer to~\cite{%
beck2020overcoming}
for MLP approximation methods
for heat PDEs with non-globally Lipschitz continuous nonlinearities.
We refer to~\cite{%
beck2020arxiv2003}
for MLP approximation methods
for semilinear elliptic PDEs.
We refer to~\cite{beck2020nested}
for MLP approximation methods for
nested conditional expectations.
We refer to~\cite{%
E2019multilevel,becker2020numerical}
for numerical simulations of MLP approximation methods
for semilinear parabolic PDEs.
We also refer to the article
E et al.~\cite{e2020algorithms} for an overview
of numerical approximation methods for high-dimensional PDEs.

\section{Simulations}
\label{sec:simulations}

In this section we give an illustration of the capabilities of the 
two approximation methods for nonlinear PDEs presented in \cref{subsec:dgm,sec:deepsplitting}
above by using them to approximate solutions to a sine-Gordon-type PDE of the form
\begin{equation}
  \label{eq:sinegordon}
  \tfrac{\partial u}{\partial t}(t,x) 
  = 
  \Delta_x u(t,x) + \sin(u(t,x))
\end{equation}
for $(t,x)\in(0,T] \times \R^d$
using three different initial values in up to $d = 1000$ dimensions.
All the simulations in this section were run on an \textsc{NVIDIA GeForce RTX 3090} GPU
with 24 GB of graphics memory, running in a containerized \textsc{Ubuntu 18.04} instance 
using 9 cores 
of an \textsc{AMD Ryzen Threadripper PRO 3975WX} with 37 GB of system memory.

We note that our implementations of both methods are rather na\"ive, meant to
provide easy-to-understand implementations for illustration purposes of these methods, 
not to obtain state-of-the-art
performance. As mentioned in \cref{sec:other} above, for the deep Galerkin
method in particular, a large volume of literature exists regarding
DNN architectures and training strategies to improve performance.
We also refer to \cref{sec:other} above for pointers to the literature containing
extensive numerical simulations exploring the capabilities of these methods.

We remark that the generally lower accuracy of the deep Galerkin method
when compared to the deep splitting method in our simulations
is to be expected: While the deep splitting method is tailored to the
type of semilinear parabolic PDEs we approximate in this section, the
deep Galerkin method uses none of the structure of the PDE problems in question.
In fact, this is part of what makes this method so appealing: In principle,
it can be applied to nearly any kind of PDE problem and it is often quite
trivial to write down a first version of an approximation algorithm, though
getting it to converge reliably to a good approximation of the solution can
be a significant challenge.

\subsection{Simulations for the deep Galerkin method}
\label{sec:simdgm}

In this section we present the numerical simulation results of
the PDE in \cref{eq:sinegordon} in dimensions $d\in\{1,2,5,10,20,50,\allowbreak 100,\allowbreak 200\}$
with three different initial values using the deep Galerkin method
sketched in \cref{subsec:dgm} above.
For every $d\in\{1,2,5,10,20,50,100,200\}$ we approximated the solution 
$
  u=(u(t,x))_{(t,x)\in[0,\nicefrac12]\times\R^d}\in C([0,\nicefrac12]\times\R^d,\R) 
$
to \cref{eq:sinegordon} at time $t=\nicefrac 12$ and at $x=0\in\R^d$
with the following initial values:
\begin{equation}
  \label{eq:initial_values}
  \begin{aligned}
    \R^d\ni x&\mapsto u(0,x)=\sqrt{1+\norm{x}^2}\in\R
    \\
    \R^d\ni x&\mapsto u(0,x)=2/(4+\norm{x}^2)\in\R
    \\
    \R^d\ni x&\mapsto u(0,x)=\arctan\bigl(\tfrac{\norm{x}}2\bigr)\in\R
  \end{aligned}
\end{equation}
To this end, we trained for every $d\in\{1,2,5,10,20,50,100,200\}$
a DNN with $3$ hidden layers with $d+1$ neurons in the input layer,
$d+40$ layers in each of the three hidden layers, and $1$ neuron in the output layer
using the method sketched in \cref{subsec:dgm}. In each of the hidden layers, the 
affine linear transformation was followed by a multidimensional version of the
Mish activation function 
\begin{equation}
  \R\ni x\mapsto
x\tanh(\ln(1+e^x)) \in\R
\end{equation}
(cf.~Misra~\cite{misra2019mish}).
The weights of the affine linear transformations were initialized using
Xavier  initialization 
(sometimes also referred to as Glorot initialization)
based on a uniform distribution 
(cf.~Glorot \& Bengio~\cite{glorot2010understanding}).
For the training, we used the Adam optimization algorithm 
(cf.~Kingma \& Ba~\cite{kingma2014adam}) with a learning rate decreasing
exponentially
over the course of the training. More precisely, with an \emph{initial learning rate}
$\lambda\in(0,\infty)$ and a \emph{learning rate decay} $\gamma\in(0,1]$ the
learning rate at step $k$ of the training was $\lambda\gamma^k$. The values
for the initial learning rate and learning rate decay, as well as the
number of training steps were chosen depending on the dimension $d$ and can be
found in \cref{table:hyper-1,table:hyper-2}. For all initial values and
all dimensions, we used minibatches of size $256$.
The collocation points in each training step 
were sampled from a normal distribution with mean $0\in\R^d$ and standard deviation
depending on the training step. More specifically, in the $k$-th of $n$ total training steps
the collocation points were sampled from a normal distribution with mean $0\in\R^d$
and covariance $2^{\frac{6(n-k)}n-1}I_{\R^{d\times d}}$.
\medskip

\begingroup
  \begin{center}
    \begin{tabular}{rrrrrll}
    \toprule
    \small \makecell{$d$}& \small\makecell{Batch\\size} & \small\makecell{Hidden\\layers} & \small\makecell{Neurons\\per hidden\\layer} & \small Steps & \small\makecell{Initial\\learning\\rate} & \small\makecell{Learning\\rate\\decay}\\
    \midrule
    $1$ & $256$ & $3$ & $41$ & $350$ & \quad\,  $0.05$ & \quad\,  $0.98$ \\
    $2$ & $256$ & $3$ & $42$ & $350$ &  \quad\, $0.05$ &  \quad\, $0.98$ \\
    $5$ & $256$ & $3$ & $45$ & $350$ &  \quad\, $0.05$ & \quad\,  $0.98$ \\
    $10$ & $256$ & $3$ & $50$ & $350$ &  \quad\, $0.05$ & \quad\,  $0.98$ \\
    $20$ & $256$ & $3$ & $60$ & $500$ &  \quad\, $0.025$ & \quad\,  $0.99$ \\
    $50$ & $256$ & $3$ & $90$ & $750$ &  \quad\, $0.01$ &  \quad\, $0.995$ \\
    $100$ & $256$ & $3$ & $140$ & $750$ &  \quad\, $0.01$ & \quad\,  $0.995$ \\
    $200$ & $256$ & $3$ & $240$ & $750$ &  \quad\, $0.005$ & \quad\,  $0.995$ \\
    \bottomrule
  \end{tabular}
  \captionof{table}{\label{table:hyper-1}The hyperparameters used in the training for the first and third
  simulation using the
  deep Galerkin method (see the results in \cref{table:dgm-1,table:dgm-3}). }
\end{center}
\endgroup

\medskip

\needspace{5cm}

\begingroup
  \begin{center}
    \begin{tabular}{rrrrrll}
    \toprule
    \small \makecell{$d$} & \small\makecell{Batch\\size} & \small\makecell{Hidden\\layers} & \small\makecell{Neurons\\per hidden\\layer} & \small Steps & \small\makecell{Initial\\learning\\rate} & \small\makecell{Learning\\rate\\decay}\\
    \midrule
    $1$ & $256$ & $3$ & $41$ & $500$ &  \quad\, $0.05$ &  \quad\, $0.99$ \\
    $2$ & $256$ & $3$ & $42$ & $500$ &  \quad\, $0.05$ &  \quad\, $0.99$ \\
    $5$ & $256$ & $3$ & $45$ & $500$ &  \quad\, $0.05$ &  \quad\, $0.99$ \\
    $10$ & $256$ & $3$ & $50$ & $750$ &  \quad\, $0.025$ & \quad\,  $0.995$ \\
    $20$ & $256$ & $3$ & $60$ & $750$ &  \quad\, $0.025$ &  \quad\, $0.995$ \\
    $50$ & $256$ & $3$ & $90$ & $1000$ &  \quad\, $0.01$ &  \quad\, $0.998$ \\
    $100$ & $256$ & $3$ & $140$ & $1000$ & \quad\,  $0.01$ &  \quad\, $0.998$ \\
    $200$ & $256$ & $3$ & $240$ & $1000$ &  \quad\, $0.005$ & \quad\,  $0.998$ \\
    \bottomrule
  \end{tabular}
  \captionof{table}{\label{table:hyper-2}The hyperparameters used in the training for the second
  simulation using the
  deep Galerkin method (see the results in \cref{table:dgm-2}).}
\end{center}
\endgroup

\medskip

For each of the three initial values in \cref{eq:initial_values} above and for each
$d\in\{1, 2,\allowbreak  5,\allowbreak  10,\allowbreak  20,\allowbreak  50, 100, 200\}$ we ran the deep Galerkin approximation algorithm 20 times
independently and
recorded the approximation of the solution $u$ of the PDE in \cref{eq:sinegordon}
with the given initial value at $t=\nicefrac12$ and $x=0\in\R^d$. The mean and standard
deviation of these values over the 20 independent runs is given in the second and third columns,
respectively, of 
\cref{table:dgm-1,table:dgm-2,table:dgm-3} below. The fourth column gives the 
reference value computed by the MLP method using the code presented in
Becker et al.~\cite{becker2020numerical}. The fifth column gives an estimation
of the absolute $L^1$-error of the approximation, computed as the mean over the 20 runs of the 
distances of the approximation from the reference value.
The estimation of the relative $L^1$-error in the sixth column is computed as 
the absolute $L^1$-error divided by the absolute value of the reference value.
Finally, the average of the runtimes in seconds of the algorithm over the 20 runs
is given in the seventh column of \cref{table:dgm-1,table:dgm-2,table:dgm-3}.
\medskip

  \begingroup
  \begin{center}
    \begin{tabular}{rrrrrrr}
    \toprule
    \small\makecell{$d$} & \small  \makecell{Mean} & \small \makecell{Standard\\deviation} & \small \makecell{Reference\\value} & \small \makecell{Absolute\\$L^1$-error} & \small \makecell{Relative\\$L^1$-error} & \small \makecell{Average\\runtime}\\
    \midrule
    $1$ & $1.788665$ & $0.013821$ & $1.836708$ & $0.048043$ & $0.026157$ & $2.48$ \\
    $2$ & $2.048171$ & $0.025434$ & $2.109239$ & $0.061302$ & $0.029064$ & $3.47$ \\
    $5$ & $2.684484$ & $0.071031$ & $2.648632$ & $0.057190$ & $0.021592$ & $6.76$ \\
    $10$ & $3.171959$ & $0.027321$ & $3.208623$ & $0.039949$ & $0.012451$ & $17.84$ \\
    $20$ & $4.041396$ & $0.112243$ & $4.107034$ & $0.088848$ & $0.021633$ & $34.38$ \\
    $50$ & $7.193619$ & $0.120224$ & $7.490964$ & $0.297344$ & $0.039694$ & $133.41$ \\
    $100$ & $9.979369$ & $0.409066$ & $9.808476$ & $0.372600$ & $0.037988$ & $292.10$ \\
    $200$ & $13.430278$ & $0.247504$ & $14.604635$ & $1.174357$ & $0.080410$ & $673.06$ \\
    \bottomrule
  \end{tabular}
  \captionof{table}{\label{table:dgm-1}Approximations for $u(\nicefrac12, 0, 0, \dots, 0)$ 
  where $u$ is the solution of the PDE 
  in \cref{eq:sinegordon}
  with the initial value $\R^d\ni x\mapsto u(0,x)=\sqrt{1+\norm{x}^2}\in\R$ for
  $d\in\{1,2,5,10,20,50,100,200\}$
  using the deep Galerkin method.
  For the hyperparameters used in the training, see \cref{table:hyper-1}.
  For the \Python\ source code used to obtain these results, see \cref{sec:source_dgm}.
  }
\end{center}
\endgroup

\medskip

\needspace{5cm}

  \begingroup
  \begin{center}
    \begin{tabular}{rrrrrrr}
    \toprule
    \small\makecell{$d$} & \small  \makecell{Mean} & \small \makecell{Standard\\deviation} & \small \makecell{Reference\\value} & \small \makecell{Absolute\\$L^1$-error} & \small \makecell{Relative\\$L^1$-error} & \small \makecell{Average\\runtime}\\
    \midrule
    $1$ & $0.680122$ & $0.003512$ & $0.677511$ & $0.003552$ & $0.005243$ & $3.62$ \\
    $2$ & $0.590964$ & $0.009817$ & $0.584510$ & $0.010066$ & $0.017221$ & $4.87$ \\
    $5$ & $0.374912$ & $0.033546$ & $0.404080$ & $0.030547$ & $0.075597$ & $9.65$ \\
    $10$ & $0.281780$ & $0.014043$ & $0.258967$ & $0.025007$ & $0.096566$ & $26.19$ \\
    $20$ & $0.139796$ & $0.013930$ & $0.147140$ & $0.013351$ & $0.090739$ & $49.88$ \\
    $50$ & $0.076818$ & $0.005612$ & $0.063228$ & $0.013589$ & $0.214926$ & $175.57$ \\
    $100$ & $0.043529$ & $0.008632$ & $0.032302$ & $0.011817$ & $0.365843$ & $384.41$ \\
    $200$ & $0.020501$ & $0.014760$ & $0.016318$ & $0.012173$ & $0.745950$ & $895.37$ \\
    \bottomrule
  \end{tabular}
  \captionof{table}{\label{table:dgm-2}Approximations for $u(\nicefrac12, 0, 0, \dots, 0)$ 
  where $u$ is the solution of the PDE 
  in \cref{eq:sinegordon}
  with the initial value $\R^d\ni x\mapsto u(0,x)=2/(4+\norm{x}^2)\in\R$ for
  $d\in\{1,2,5,10,20,50,100,200\}$
  using the deep Galerkin method.
  For the hyperparameters used in the training, see \cref{table:hyper-2}.
  For the \Python\ source code used to obtain these results, see \cref{sec:source_dgm}.
  }
\end{center}
\endgroup

\medskip

\needspace{5cm}

  \begingroup
  \begin{center}
    \begin{tabular}{rrrrrrr}
    \toprule
    \small\makecell{$d$} & \small  \makecell{Mean} & \small \makecell{Standard\\deviation} & \small \makecell{Reference\\value} & \small \makecell{Absolute\\$L^1$-error} & \small \makecell{Relative\\$L^1$-error} & \small \makecell{Average\\runtime}\\
    \midrule
    $1$ & $0.629165$ & $0.029052$ & $0.569925$ & $0.061155$ & $0.107303$ & $2.47$ \\
    $2$ & $0.900790$ & $0.016736$ & $0.835350$ & $0.065440$ & $0.078338$ & $3.52$ \\
    $5$ & $1.249522$ & $0.020660$ & $1.201798$ & $0.047724$ & $0.039711$ & $6.77$ \\
    $10$ & $1.395483$ & $0.022569$ & $1.440293$ & $0.044975$ & $0.031226$ & $18.29$ \\
    $20$ & $1.593871$ & $0.029176$ & $1.622937$ & $0.035450$ & $0.021843$ & $34.28$ \\
    $50$ & $1.789839$ & $0.032926$ & $1.785912$ & $0.027953$ & $0.015652$ & $131.69$ \\
    $100$ & $1.851764$ & $0.034147$ & $1.866269$ & $0.030308$ & $0.016240$ & $291.67$ \\
    $200$ & $1.626139$ & $0.206072$ & $1.921878$ & $0.295739$ & $0.153880$ & $678.89$ \\
    \bottomrule
  \end{tabular}
  \captionof{table}{\label{table:dgm-3}Approximations for $u(\nicefrac12, 0, 0, \dots, 0)$ 
  where $u$ is the solution of the PDE 
  in \cref{eq:sinegordon}
  with the initial value
  $\R^d\ni x\mapsto u(0,x)=\arctan\bigl(\tfrac{\norm{x}}2\bigr)\in\R$ for
  $d\in\{1,2,5,10,20,50,100,200\}$
  using the deep Galerkin method.
  For the hyperparameters used in the training of the DNNs, see \cref{table:hyper-1}.
  For the \Python\ source code used to obtain these results, see \cref{sec:source_dgm}.
  }
\end{center}
\endgroup

\subsection{Simulations for the deep splitting method}
\label{sec:simds}
In this section we present the numerical simulation results for
the PDE in \cref{eq:sinegordon} in dimensions $d\in\{1,2,5,10,\allowbreak 20,\allowbreak 50,\allowbreak 100,200, 500, 1000\}$
with three different initial values using the deep splitting method sketched in
\cref{sec:deepsplitting} above.
For every $d\in\{1,2,5,10,\allowbreak 20,\allowbreak 50,\allowbreak 100,200, 500, 1000\}$ we
approximated the solution 
$u=(u(t,x))_{(t,x)\in[0,\nicefrac12]\times\R^d}\in C([0,\nicefrac12]\times\R^d,\R)$
to \cref{eq:sinegordon} at time $t=\nicefrac 12$ and at $x=0\in\R^d$.
To this end, we split up the time interval $[0,\nicefrac12]$ into
$30$ subintervals and successively trained $30$ artificial neural networks to
approximate the solution $u$ at each of the time 
points $\nicefrac{i}{60}$ for $i\in\{1, 2, \dots,30\}$.
Each of these artificial neural networks was a fully-connected
feedforward neural network with $2$ hidden layers with $d$ neurons in the input layer,
$d+50$ neurons in each of the two hidden layers,
and $1$ neuron in the output layer (except for the one used to approximate the 
solution at time $T=\nicefrac12$, which consisted of a single affine linear transformation). 
In each of the hidden layers the affine linear
transformation was followed by batch normalization (cf.~Ioffe \& Szegedy~\cite{ioffe2015batch}) and a multidimensional
version of the exponential linear unit 
(ELU) activation function 
\begin{equation}
  \R\ni x\mapsto x + (e^x - 1 -x)\ind{(-\infty,0]}(x)\in\R
\end{equation}
(cf.\ Clevert et al.~\cite{clevert2015fast})
The weights of the affine linear transformations were initialized using
Xavier initialization 
based on a uniform distribution 
(cf.~Glorot \& Bengio~\cite{glorot2010understanding}).
As for the simulations described in \cref{sec:simdgm}, 
we used the Adam optimization algorithm with a learning rate decreasing exponentially
over the course of the training. The initial learning rate was
$0.2$ for all simluations.
The learning rate decay, as well as the
number of training steps were chosen depending on the dimension $d$ and can be
found in \cref{table:hyper-1,table:hyper-2}. For all initial values and
all dimensions, we used minibatches of size $256$.
The sampling of the collocation points used during training
was implemented as described in \cite{beck2021deep}. More specifically,
to train the $k$-th neural network, i.e., the network used to approximate the
solution of the PDE in question at time $\frac{kT}n$,
the collocation points were sampled from a normal distribution with mean $0\in\R^d$ and
covariance $(1-\frac k{30})I_{\R^{d\times d}}$.

\medskip

\begingroup
\begin{minipage}{\linewidth}
  \begin{center}
  \begin{tabular}{rrrrrll}
    \toprule
    \small $d$ & \small\makecell{Batch\\size} & \small\makecell{Hidden\\layers} & \small\makecell{Neurons\\per hidden\\layer} & \small Steps & \small\makecell{Initial\\learning\\rate} & \small\makecell{Learning\\rate\\decay}\\
    \midrule
    $1$ & $256$ & $2$ & $51$ & $350$ & \quad\,  $0.2$ & \quad\,  $0.985$ \\
    $2$ & $256$ & $2$ & $52$ & $350$ &  \quad\, $0.2$ &  \quad\, $0.985$ \\
    $5$ & $256$ & $2$ & $55$ & $350$ &  \quad\, $0.2$ & \quad\,  $0.985$ \\
    $10$ & $256$ & $2$ & $60$ & $350$ &  \quad\, $0.2$ & \quad\,  $0.985$ \\
    $20$ & $256$ & $2$ & $70$ & $500$ &  \quad\, $0.2$ & \quad\,  $0.99$ \\
    $50$ & $256$ & $2$ & $100$ & $500$ &  \quad\, $0.2$ &  \quad\, $0.99$ \\
    $100$ & $256$ & $2$ & $150$ & $750$ &  \quad\, $0.2$ & \quad\,  $0.995$ \\
    $200$ & $256$ & $2$ & $250$ & $1000$ &  \quad\, $0.2$ & \quad\,  $0.995$ \\
    $500$ & $256$ & $2$ & $550$ & $1000$ &  \quad\, $0.2$ & \quad\,  $0.995$ \\
    $1000$ & $256$ & $2$ & $1050$ & $1250$ &  \quad\, $0.2$ & \quad\,  $0.998$ \\
    \bottomrule
  \end{tabular}
  \captionof{table}{\label{table:hyper-ds}The hyperparameters used in the training for the 
  simulations using the
  deep splitting method (see the results in \cref{table:ds-1,table:ds-2,table:ds-3}) }
\end{center}
\end{minipage}
\endgroup

\medskip

For each of the three initial values in \cref{eq:initial_values} above and for each
$d\in\{1, 2, 5, 10, 20,\allowbreak 50,\allowbreak 100, 200, 500, 1000\}$ we ran the deep splitting 
approximation algorithm 20 times independently and
recorded the approximation of the solution $u$ of the PDE in \cref{eq:sinegordon}
with the given initial value at $t=\nicefrac12$ and $x=0\in\R^d$. 
The columns of the results tables, \cref{table:ds-1,table:ds-2,table:ds-3} below,
follow the same format as in \cref{sec:simdgm}.

\medskip

\needspace{2cm}

\begingroup
\begin{minipage}{\linewidth}
  \begin{center}
    \begin{tabular}{rrrrrrr}
    \toprule
    \small\makecell{$d$} & \small  \makecell{Mean} & \small \makecell{Standard\\deviation} & \small \makecell{Reference\\value} & \small \makecell{Absolute\\$L^1$-error} & \small \makecell{Relative\\$L^1$-error} & \small \makecell{Average\\runtime}\\
    \midrule
    $1$ & $1.830958$ & $0.015773$ & $1.836708$ & $0.013452$ & $0.007324$ & $20.22$ \\
    $2$ & $2.108470$ & $0.014779$ & $2.109239$ & $0.012182$ & $0.005776$ & $19.90$ \\
    $5$ & $2.642265$ & $0.015791$ & $2.648632$ & $0.014881$ & $0.005618$ & $20.13$ \\
    $10$ & $3.196520$ & $0.016279$ & $3.208623$ & $0.016151$ & $0.005034$ & $20.31$ \\
    $20$ & $4.082424$ & $0.019874$ & $4.107034$ & $0.026197$ & $0.006379$ & $28.89$ \\
    $50$ & $7.480304$ & $0.027910$ & $7.490964$ & $0.022157$ & $0.002958$ & $29.73$ \\
    $100$ & $9.804550$ & $0.009534$ & $9.808476$ & $0.008295$ & $0.000846$ & $46.51$ \\
    $200$ & $14.602816$ & $0.036716$ & $14.604635$ & $0.020845$ & $0.001427$ & $64.86$ \\
    $500$ & $21.386109$ & $1.045499$ & $22.230142$ & $0.844033$ & $0.037968$ & $75.23$ \\
    $1000$ & $31.748045$ & $0.126012$ & $31.751259$ & $0.068201$ & $0.002148$ & $117.77$ \\
    \bottomrule
  \end{tabular}
  \captionof{table}{\label{table:ds-1}Approximations for $u(\nicefrac12, 0, 0, \dots, 0)$ 
  where $u$ is the solution of the PDE 
  in \cref{eq:sinegordon}
  with the initial value $\R^d\ni x\mapsto u(0,x)=\sqrt{1+\norm{x}^2}\in\R$ for
  $d\in\{1,2,5,10,20,50,100,200,500,1000\}$
  using the deep splitting method.
  For the hyperparameters used in the training of the DNNs, see \cref{table:hyper-ds}.
  For the \Python\ source code used to obtain these results, see \cref{sec:source_ds}.
  }
\end{center}
\end{minipage}
\endgroup

\medskip

  \begingroup
  \begin{center}
  \begin{tabular}{rrrrrrr}
    \toprule
    \small\makecell{$d$} & \small  \makecell{Mean} & \small \makecell{Standard\\deviation} & \small \makecell{Reference\\value} & \small \makecell{Absolute\\$L^1$-error} & \small \makecell{Relative\\$L^1$-error} & \small \makecell{Average\\runtime}\\
    \midrule
    $1$ & $0.678383$ & $0.004441$ & $0.677511$ & $0.003652$ & $0.005390$ & $20.24$ \\
    $2$ & $0.583808$ & $0.005164$ & $0.584510$ & $0.004229$ & $0.007235$ & $19.93$ \\
    $5$ & $0.405637$ & $0.003722$ & $0.404080$ & $0.003680$ & $0.009107$ & $20.08$ \\
    $10$ & $0.257737$ & $0.001830$ & $0.258967$ & $0.001949$ & $0.007524$ & $20.33$ \\
    $20$ & $0.146415$ & $0.000674$ & $0.147140$ & $0.000859$ & $0.005841$ & $29.08$ \\
    $50$ & $0.062905$ & $0.000134$ & $0.063228$ & $0.000323$ & $0.005114$ & $30.05$ \\
    $100$ & $0.032181$ & $0.000150$ & $0.032302$ & $0.000146$ & $0.004514$ & $46.48$ \\
    $200$ & $0.016254$ & $0.000023$ & $0.016318$ & $0.000064$ & $0.003941$ & $64.93$ \\
    $500$ & $0.006542$ & $0.000005$ & $0.006568$ & $0.000026$ & $0.003909$ & $75.37$ \\
    $1000$ & $0.003277$ & $0.000010$ & $0.003291$ & $0.000015$ & $0.004467$ & $118.31$ \\
    \bottomrule
  \end{tabular}
  \captionof{table}{\label{table:ds-2}Approximations for $u(\nicefrac12, 0, 0, \dots, 0)$ 
  where $u$ is the solution of the PDE 
  in \cref{eq:sinegordon}
  with the initial value $\R^d\ni x\mapsto u(0,x)=2/(4+\norm{x}^2)\in\R$ for
  $d\in\{1,2,5,10,20,50,100,200,500,1000\}$
  using the deep splitting method.
  For the hyperparameters used in the training of the DNNs, see \cref{table:hyper-ds}.
  For the \Python\ source code used to obtain these results, see \cref{sec:source_ds}.
  }
\end{center}
\endgroup

\medskip

  \begingroup
  \begin{center}
  \begin{tabular}{rrrrrrr}
    \toprule
    \small\makecell{$d$} & \small  \makecell{Mean} & \small \makecell{Standard\\deviation} & \small \makecell{Reference\\value} & \small \makecell{Absolute\\$L^1$-error} & \small \makecell{Relative\\$L^1$-error} & \small \makecell{Average\\runtime}\\
    \midrule
    $1$ & $0.561503$ & $0.012954$ & $0.569925$ & $0.012631$ & $0.022162$ & $20.16$ \\
    $2$ & $0.830821$ & $0.011572$ & $0.835350$ & $0.009811$ & $0.011745$ & $20.19$ \\
    $5$ & $1.195350$ & $0.007251$ & $1.201798$ & $0.008002$ & $0.006658$ & $20.09$ \\
    $10$ & $1.436600$ & $0.003500$ & $1.440293$ & $0.004196$ & $0.002913$ & $20.22$ \\
    $20$ & $1.622351$ & $0.001493$ & $1.622937$ & $0.001305$ & $0.000804$ & $29.17$ \\
    $50$ & $1.785940$ & $0.000271$ & $1.785912$ & $0.000214$ & $0.000120$ & $29.81$ \\
    $100$ & $1.866438$ & $0.000438$ & $1.866269$ & $0.000337$ & $0.000181$ & $46.65$ \\
    $200$ & $1.922286$ & $0.000091$ & $1.921878$ & $0.000408$ & $0.000212$ & $65.01$ \\
    $500$ & $1.970859$ & $0.000031$ & $1.970294$ & $0.000566$ & $0.000287$ & $75.16$ \\
    $1000$ & $1.994976$ & $0.000086$ & $1.994330$ & $0.000645$ & $0.000324$ & $118.21$ \\
    \bottomrule
  \end{tabular}
  \captionof{table}{\label{table:ds-3}Approximations for $u(\nicefrac12, 0, 0, \dots, 0)$ 
  where $u$ is the solution of the PDE 
  in \cref{eq:sinegordon}
  with the initial value $\R^d\ni x\mapsto u(0,x)=\arctan\bigl(\tfrac{\norm{x}}2\bigr)\in\R$ for
  $d\in\{1,2,5,10,20,50,100,200,500,1000\}$
  using the deep splitting method.
  For the hyperparameters used in the training of the DNNs, see \cref{table:hyper-ds}.
  For the \Python\ source code used to obtain these results, see \cref{sec:source_ds}.
  }
\end{center}
\endgroup

\section{Theoretical results for DNN approximations for PDEs}
\label{sec:theoretical}

In this section we will briefly mention some of the results that exist 
in the scientific literature which attempt to explain in theoretical terms
the success of deep learning-based approximation methods 
for high-dimensional PDEs in numerical simulations.
We first provide a general introduction to the topic (see \cref{sec:theorylit})
and then review the main result in Hutzenthaler et al.~\cite{hutzenthaler2020proof} 
as one instance of a mathematical result providing a partial error bound 
for deep learning-based approximation methods for PDEs (see \cref{sec:approxerror}).

\subsection[Literature overview of theoretical results for DNN approximations for PDEs]{Literature overview of theoretical results for DNN approximations for high-dimensional PDEs}
\label{sec:theorylit}

Although by now a large number of highly encouraging
numerical simulations hinting at the great potential
of deep learning-based methods in the approximation
of high-dimensional PDEs can be found in the scientific
literature (see the references in \cref{sec:other}; 
cf.~\cref{sec:simulations} above), 
there is as yet no mathematical result
rigorously explaining the great practical success 
of deep learning-based approximation methods for PDEs.
In particular, so far there does not exist a proof
that deep learning-based methods
are capable of overcoming the curse of dimensionality
in the approximative computation of solutions of
high-dimensional PDEs in the sense that the number of
computational operations necessary for calculating the numerical 
approximation of the PDE under consideration grows
at most polynomially in the PDE dimension $d\in\N$
and the reciprocal $\nicefrac1\eps$ of the
prescribed approximation accuracy $\eps\in(0,\infty)$.

However, recent results supplying
rigorous partial error analyses show that DNNs
at least possess the necessary expressive power to overcome
the curse of dimensionality in the approximation of
high-dimensional PDEs in the sense that there
exist DNN approximations for 
high-dimensional PDEs with the number of
parameters used to describe the approximating
neural network growing at most polynomially in the
PDE dimension $d\in\N$ and the reciprocal $\nicefrac1\eps$
of the prescribed approximation accuracy $\eps\in(0,\infty)$.
The first result in this direction was obtained in
2018 in Grohs et al.~\cite{grohs2018proof}
and, since then, a number of articles have appeared
significantly extending the results in Grohs et al.~\cite{grohs2018proof},
covering several more classes of PDEs (see, e.g., \cite{BernerGrohsJentzen2018,
elbraechter2018dnn,
gonon2019uniform,
GrohsHerrmann2020,
grohs2019deep,
HornungJentzenSalimova2020,
jentzen2021proof,
reisinger2019rectified,
hutzenthaler2020proof,
gonon2020deep,
kutyniok2019theoretical,
beneventano2020highdimensional}).
In particular, in \cite{hutzenthaler2020proof} there is now also a result  
on nonlinear PDEs in the scientific literature
which shows that DNNs have the expressive power to overcome
the curse of dimensionality in the numerical approximation of
solutions of nonlinear heat PDEs with Lipschitz continuous nonlinearities.

To give the reader a somewhat more complete picture on theoretical
results for DNN approximations for PDEs, we also refer, e.g., to
Han \& Long \cite{HanLong2020} and Hur\'e et al.~\cite{HurePhamWarin2019} for 
conditional convergence rate results
for approximation errors
for deep learning-based approximations for PDEs,
we also refer, e.g., to Luo \& Yang~\cite{LuoYang2020},
Shin et al.~\cite{shin2020convergence},
Mishra \& Molinaro~\cite{mishra2021estimates},
and Hu et al.~\cite{hu2022when}
 for upper bounds
for the generalization error for deep learning-based approximations
for PDEs,
we also refer, e.g., to Darbon et al.~\cite{DarbonLangloisMeng2020}
and Darbon \& Meng~\cite{DarbonMeng2021}
for exact DNN representation results for certain 
Hamilton--Jacobi PDEs,
and we also refer, e.g., to Geist et al.~\cite{GeistEtal2020} for numerical
investigations regarding the expressive power of
DNNs in the context of the approximation of 
solutions of PDEs.

There are to our best knowledge no results to date rigorously showing
for a neural network-based approximation method for PDEs 
that the optimization error associated to an SGD type method with one random initialization
is not subject to the curse of
dimensionality, nor that it even converges to zero, except for the special
situation where only the parameters
in the output layer are learned during training, in which case
it has been shown under suitable assumptions that
shallow neural networks can learn the solutions to
Black--Scholes type PDEs without the curse of dimensionality
(see Gonon~\cite{gonon2021random}).

\subsection{Approximation error bounds for DNN approximations
for nonlinear heat PDEs}
\label{sec:approxerror}

In order to provide the reader with a more formal illustration of the type of 
results obtained in the branch of research outlined in
\cref{sec:theorylit} above, we present in this section
in the following theorem, \cref{thm:DNN_nonlinear_PDEs} below,
a special case of
the main result in \cite{hutzenthaler2020proof}.
\cfclear
\begin{theorem}
  \label{thm:DNN_nonlinear_PDEs}
  Let $
    {\bf N} = \allowbreak
    \bigcup_{ L \in \N }\allowbreak
    \bigcup_{ l_0, l_1, \dots, l_L \in \N }(
      \bigtimes_{ k = 1 }^L ( \R^{ l_k \times l_{k-1} } \times \R^{l_k} )
    ) 
  $, 
  let 
  $ \mathcal{R} \colon {\bf N} \to (\bigcup_{k,l\in\N}C(\R^k,\R^l)) $ 
  and 
  $
    \mathcal{P} \colon {\bf N} \to \N 
  $ 
  satisfy for all 
    $L\in\N$, 
    $l_0,l_1,\dots,l_L\in\N$, 
    $\Phi = ((W_1,B_1),\allowbreak(W_2,B_2),\allowbreak\dots,\allowbreak(W_L,B_L)) \in (\bigtimes_{k=1}^L\allowbreak (\R^{l_k\times l_{k-1}}\times \R^{l_k}))$, 
    $x_0\in \R^{l_0}$, $x_1\in\R^{l_1}$, $\dots$, $x_L\in\R^{l_L}$ with 
    $\forall\, k\in \N\cap(0,L)\colon x_k = \rect{l_k}( W_k x_{ k - 1 } + B_k ) $ that
  $
    \mathcal{R}(\Phi) \in C(\R^{l_0},\R^{l_{L}})
  $, 
  $
    (\mathcal{R}(\Phi))(x_0) = W_Lx_{L-1}+B_L
  $, 
  and
  $
    \mathcal{P}(\Phi) = \sum_{k=1}^L l_k(l_{k-1}+1)
  $, 
  let $ T, \kappa \in (0,\infty) $, 
  $ (\mathfrak{g}_{d,\varepsilon} )_{ (d,\varepsilon) \in \N \times (0,1] } \subseteq {\bf N} $, 
  let $u_d \in C^{1,2}([0,T]\times \R^d,\R)$, $d\in\N$,
  let $f\colon\R \to \R$ be Lipschitz continuous, 
  and assume for all 
  $ d \in \N $, $ x = (x_1, x_2, \dots, x_d) \in \R^d $, 
  $ \varepsilon \in (0,1] $, $ t \in [0,T] $ 
  that 
  $
    \mathcal{R}(\mathfrak{g}_{d,\varepsilon}) \in C(\R^d,\R) 
  $, 
  $ 
    \varepsilon \abs{ u_d(t,x) }
    + \allowbreak 
    \abs{ u_d(0,x) \allowbreak - ( \mathcal{R}(\mathfrak{g}_{d,\varepsilon}) )(x) }
    \le \varepsilon \kappa d^\kappa (1 + \sum_{ i = 1 }^d \abs{ x_i }^\kappa ) 
  $, 
  $
    \mathcal{P}(\mathfrak{g}_{d,\varepsilon}) \le \kappa d^\kappa \varepsilon^{-\kappa} 
  $, 
  and 
  \begin{equation}\label{eq:1a}
    \tfrac{\partial u_d}{\partial t}(t,x) = \Delta_x u_d(t,x) + f(u_d(t,x))
  \end{equation}
  \cfload.
  Then there exist 
  $ (\mathfrak{u}_{d,\varepsilon})_{(d,\varepsilon)\in \N\times(0,1] }\subseteq {\bf N} $
  and $c \in \R $ 
  such that for all $ d \in \N $, $ \varepsilon \in (0,1] $ 
  it holds that 
  $ \mathcal{R}( \mathfrak{u}_{d,\varepsilon} ) \in C( \R^d, \R ) $, 
  $ \mathcal{P}( \mathfrak{u}_{d,\varepsilon} ) \le c d^c \varepsilon^{-c} $, 
  and
  \begin{equation}\label{eq:2a}
    \br*{
      \int_{ [0,1]^d } 
        \abs*{
          u_d(T,x) - (\mathcal{R}(\mathfrak{u}_{d,\varepsilon}))(x) 
        }^2 
      \,dx
    }^{\nicefrac{1}{2}} \le \varepsilon.
  \end{equation}
\end{theorem}
\cref{thm:DNN_nonlinear_PDEs} follows immediately from~\cite[Theorem~1.1]{hutzenthaler2020proof}.
In the following we add some comments on the mathematical objects appearing in
\cref{thm:DNN_nonlinear_PDEs} above.
The set $\mathbf N$ in \cref{thm:DNN_nonlinear_PDEs} represents the set
of all neural networks. More specifically, each element 
$\Phi\in\pr[\big]{\bigtimes_{k=1}^L (\R^{l_k\times l_{k-1}}\times \R^{l_k})}$
with $L\in\N$, $l_0,l_1,\dots,l_L\in\N$ represents a neural network
where $L\in\N$ is the length of the neural network and $l_0,l_1,\dots,l_L\in\N$
correspond to the number of neurons in the $1$st, $2$nd, \ldots, $(L+1)$th layer,
respectively.
The function 
\begin{equation}
  \textstyle
\mathcal{R} \colon {\bf N} \to \bigl(\bigcup_{k,l\in\N}C(\R^k,\R^l)\bigr)
\end{equation}
in \cref{thm:DNN_nonlinear_PDEs} assigns to each neural network its
associated realization function.
More specifically, for each neural network $\Phi\in\mathbf N$
it holds that 
$\mathcal R(\Phi)\in (\bigcup_{k,l\in\N}C(\R^k,\R^l))$ is the realization
function of the neural network $\Phi$, where the activation functions
are multidimensional versions of the
rectifier function 
introduced in \cref{def:rect} above.
The function $\mathcal P\colon \mathbf N\to\N$ assigns to each
neural network $\Phi\in\mathbf N$ the number of parameters employed in the description
of the neural network $\Phi$. 
We observe that for every
$\Phi\in\mathbf N$ it holds that
$\mathcal P(\Phi)$ corresponds, roughly speaking, to the amount of memory
necessary to store the neural network $\Phi$ on a computer.

In \eqref{eq:1a} in \cref{thm:DNN_nonlinear_PDEs} we specify the
PDEs whose solutions we intend to approximate by neural networks.
The functions $u_d\colon[0,T]\times\R^d\to\R$, $d\in\N$,
in \cref{thm:DNN_nonlinear_PDEs} denote the exact solutions
of the PDEs in \eqref{eq:1a}.
The real number $T\in(0,\infty)$ in \cref{thm:DNN_nonlinear_PDEs}
denotes the time horizon of the
PDEs in \eqref{eq:1a}.
The function $f\colon\R\to\R$ in
\cref{thm:DNN_nonlinear_PDEs}
describes the Lipschitz nonlinearity
in the PDEs in \eqref{eq:1a}.
The real number $\kappa\in(0,\infty)$ in
\cref{thm:DNN_nonlinear_PDEs} is a constant used to
formulate the regularity and approximation hypotheses
in \cref{thm:DNN_nonlinear_PDEs}.
We assume in \cref{thm:DNN_nonlinear_PDEs} that the solutions 
$u_d\colon[0,T]\times\R^d\to\R$, 
$d\in\N$, of the PDEs in \eqref{eq:1a}
grow at most polynomially.
The real number $\kappa\in(0,\infty)$
is used to formulate this polynomial growth
assumption in \cref{thm:DNN_nonlinear_PDEs}.
More formally, in \cref{thm:DNN_nonlinear_PDEs}
it is assumed that for all $d\in\N$,
$t\in[0,T]$,
$x=(x_1,x_2,\dots,x_d)\in\R^d$ 
it holds that 
\begin{equation}
  \textstyle
\abs{u_d(t,x)}
\leq 
\kappa d^\kappa\bigl(1+\sum_{k=1}^d\abs{x_i}^\kappa\bigr).
\end{equation}
The neural networks $\mathfrak g_{d,\eps}\in\mathbf N$, $d\in\N$, $\eps\in(0,1]$,
in \cref{thm:DNN_nonlinear_PDEs}
describe neural network approximations to the initial conditions
$\R^d\ni x\mapsto u_d(0,x)\in \R$, $d\in\N$, of the PDEs
in \eqref{eq:1a}.
In particular, note that the hypothesis in
\cref{thm:DNN_nonlinear_PDEs} that for all
$d\in\N$, $x=(x_1,x_2,\dots,x_d)\in\R^d$,
$\eps\in(0,1]$, $t\in[0,T]$ it holds that
\begin{equation}
  \eps\abs{u_d(t,x)}+\abs{u_d(0,x)-(\mathcal R(\mathfrak g_{d,\eps}))(x)}
  \leq 
  \eps\kappa d^\kappa\pr[\big]{1+{\textstyle\sum}_{i=1}^d\abs{x_i}^\kappa}
\end{equation}
implies that for all $d\in\N$, $x\in\R^d$ it holds that
$(\mathcal R(\mathfrak g_{d,\eps}))(x)$ converges to $u_d(0,x)$ 
as $\eps$ tends to $0$.
Observe that this in combination with the assumption in \cref{thm:DNN_nonlinear_PDEs}
that for all $d\in\N$, $\eps\in(0,1]$ it holds that
$\mathcal P(\mathfrak g_{d,\eps})\leq \kappa d^\kappa\eps^{-\kappa}$ ensures
that the initial conditions of the PDEs in \eqref{eq:1a}
can be approximated through neural networks without the curse
of dimensionality.

\cref{thm:DNN_nonlinear_PDEs} establishes that there exist neural networks
$\mathfrak u_{d,\eps}\in\mathbf N$, $d\in\N$, $\eps\in(0,1]$, 
such that for all $d\in\N$, $\eps\in(0,1]$ it holds that the $L^2$-distance,
with respect to the Lebesgue measure on the region $[0,1]^d$,
between the exact solution $\R^d\ni x\mapsto u_d(T,x)\in\R$ of the PDE 
in \eqref{eq:1a} at the terminal time
$T\in (0,\infty)$ and the realization function 
$\mathcal R(\mathfrak u_{d,\eps})\colon \R^d\to\R$ of the neural network
$\mathfrak u_{d,\eps}\in\mathbf N$
is bounded by $\eps$ and such that the number of parameters of the neural networks 
$\mathfrak u_{d,\eps}\in\mathbf N$, $d\in\N$, $\eps\in(0,1]$, grows at most polynomially
in both the PDE dimension $d\in\N$ and the reciprocal $\nicefrac 1\eps$ of the
prescribed approximation accuracy $\eps\in(0,1]$.
More precisely, \cref{thm:DNN_nonlinear_PDEs} establishes that there
exist real numbers $K,c\in\R$ such that
it holds for all $d\in\N$, $\eps\in(0,1]$ that 
\begin{equation}
\mathcal{P}( \mathfrak{u}_{d,\varepsilon} ) \le K d^c \varepsilon^{-c}.
\end{equation}
We remark that while the statement of
\cref{thm:DNN_nonlinear_PDEs} in principle allows both
$K$ and $c$ in the claim to depend on the
time horizon $T\in(0,\infty)$, we expect that $c$ can in fact
be chosen independently of $T$.
In contrast, we expect $K$ to grow exponentially with $T$
(in the absence of additional assumptions),
similar to approximations of ODEs by the Euler method
(see also \cite[Theorem~3.8]{hutzenthaler2020overcoming}). 
It is conceivable that this
dependency could be improved, as in the approximation of
ODEs, with additional assumptions, such as dissipativity.

In the following we add some comments on the arguments of the proof of
\cref{thm:DNN_nonlinear_PDEs}. The key idea of the proof of \cref{thm:DNN_nonlinear_PDEs} is
\begin{enumerate}[(I)]
  \item[\mylabel{it:dnn1}{($\mathfrak A$)}]
  to consider, for every $\eps\in(0,1]$, perturbed versions of the PDEs in \eqref{eq:1a}
  in which the nonlinearity $f\colon\R\to\R$ is replaced
  by the realization of a suitable DNN and
  in which the initial condition functions
  $\R^d\ni x\mapsto u_d(0,x)\in\R$, $d\in\N$,
  are replaced by the realizations
  $\mathcal R(\mathfrak g_{d,\eps})\colon \R^d\to\R$, $d\in\N$,
  of the given DNNs $\mathfrak g_{d,\eps}\in\mathbf N$, $d\in\N$,
  \item[\mylabel{it:dnn2}{($\mathfrak B$)}] 
  to approximate the perturbed versions
  of the PDEs in \eqref{eq:1a}
  according to item \ref{it:dnn1}
  by means of MLP approximation methods
  (which are known to overcome the
  curse of dimensionality in the approximation of nonlinear heat
  equations with Lipschitz nonlinearities; see, 
  e.g., \cite{hutzenthaler2020overcoming} and \cref{sec:intro} above),
  and
  \item[\mylabel{it:dnn3}{($\mathfrak C$)}]
  to represent suitable random realizations
  of these MLP approximations as realization functions of DNNs by mimicking the construction
  of MLP approximations through DNNs.
\end{enumerate}
We refer to \cite[Section 7]{e2020algorithms} for further details and to
\cite{hutzenthaler2020proof} for the detailed proof
of \cref{thm:DNN_nonlinear_PDEs}.

As indicated in item \ref{it:dnn2} above, the proof of \cref{thm:DNN_nonlinear_PDEs} above crucially 
uses the fact that MLP approximation algorithms
do overcome the curse of dimensionality in the approximative
solution of nonlinear heat PDEs with Lipschitz continuous nonlinearities
(see \cite{hutzenthaler2020overcoming}).
We would like to point out that, while \cref{thm:DNN_nonlinear_PDEs} above
is a statement about the class of 
nonlinear heat PDEs with Lipschitz continuous nonlinearities,
the fact that MLP approximations overcome the curse of dimensionality
has been established for a much larger class of PDEs 
(see \cite{beck2020overcoming} for MLP approximations for
heat PDEs with possibly non-Lipschitz 
continuous nonlinearities; see
\cite{beck2020arxiv2003}
for MLP approximations for semilinear elliptic PDEs; see
\cite{hutzenthaler2019arxiv1903,hutzenthaler2020lipschitz} for MLP approximations
for parabolic PDEs with more general second order differential
operators than just the Laplacian; 
see \cite{hutzenthaler2020multilevel}
for MLP approximations for semilinear PDEs with gradient-dependent nonlinearities).
These results make it seem feasible that
results analogous to \cref{thm:DNN_nonlinear_PDEs} above may
be established for other classes of PDEs as well.

\section{Conclusion}

In recent years, deep learning has become a powerful tool for
approximating solutions of PDEs, especially in high dimensions.
In this article we have focused on describing a selection 
of deep learning-based approximation
methods for high-dimensional PDEs from the literature. We also briefly 
reviewed the wider field as well as selected other
Monte Carlo methods for approximating solutions
of high-dimensional PDEs.

We believe that these developments pave the way to an exciting future
for this area of research.
The tools surveyed here open the door to finding solutions for a
range of practically relevant problems that were inaccessible 
until quite recently. Notably, it has suddenly become feasible to 
find approximate solutions to semilinear PDEs in 1000 or 10\,000 
space dimensions 
(cf., e.g., \cite{becker2020numerical,beck2021deep} and \cref{sec:simulations} above),
a problem that until a few years ago had to be considered
completely out of reach.
In our own subjective view,
this could, on the one hand,  lead to dramatic new approaches
in applications, e.g.,
concerning
optimal control problems,
concerning nonlinear filtering problems,
concerning the approximative pricing of
products,
and concerning the modelling
of complex systems in physics and biology (see \cref{sec:intro}).
On the other hand, the developments laid out here
open up new paths of research. Firstly,
the advances reviewed here will likely
spawn new approximation methods for high-dimensional PDEs,
which will, secondly, demand rigorous mathematical
analyses of their approximative capabilities. In particular,
though the results surveyed in \cref{sec:theoretical} provide
a start to the mathematical analysis of 
deep learning-based approximation methods
for PDEs, there is still a long way to go towards
a full mathematical theory.
And thirdly, the results surveyed here might have ramifications
for a regularity theory for high-dimensional PDEs, where
higher regularity does not correspond to higher smoothness properties
of the involved functions but where higher regularity
corresponds to better approximability properties
of the involved functions in high dimensions
(cf.\ Beneventano et al.~\cite[Section~3]{beneventano2020highdimensional}).

\section{Source code for simulations}\phantom{h}
 In this appendix we provide the complete source code used to obtain the simulation
 results in \cref{sec:simulations} above.
 All of the \Python\ source files listed here can be downloaded as part of the 
 sources of the arXiv version
 of this article at \url{https://arxiv.org/e-print/2012.12348} (in the form of a gzipped
 tar file).

\subsection{Common source code}\phantom{h}

\begin{lstlisting}
import time
import json
import math 
import logging

from collections import namedtuple

JobConfig = namedtuple('JobConfig', ['build_model', 'runs', 'train', 
                                     'train_params', 'evaluate'])

JobResult = namedtuple('JobResult', ['model', 'values', 'train_times'])

def run_job(job_config):
    values = []
    train_times = []

    model = job_config.build_model()

    for run in range(job_config.runs):
        logging.debug(f"Starting Run {run+1}")

        tick = time.perf_counter()
        job_config.train(model, job_config.train_params)
        tock = time.perf_counter()

        value = job_config.evaluate(model)

        values.append(value)
        train_times.append(tock - tick)

        logging.debug(f"Value: {value}")
        logging.debug(f"Time: {tock - tick}\n")

    return JobResult(model, values, train_times)


def run_jobs(job_configs, init, report):
    results = []

    for config in job_configs:
        init()

        for (k,v) in config.train_params.items():
            logging.info(f"{k}: {v}")
        logging.info("")
        
        result = run_job(config)
        results.append(report(result, config.train_params))

    return results


def report(job_result, train_params):
    mean = sum(job_result.values) / len(job_result.values)
    result = {
        "values": job_result.values,
        "train_times": job_result.train_times,
        "model": f"{job_result.model}",
        "mean": mean, 
        "stddev": math.sqrt(sum([(y - mean)**2 for y in job_result.values]) 
                            / len(job_result.values)),
        "avgruntime": sum(job_result.train_times) / len(job_result.train_times),
        "train_params": {k:v for (k,v) in train_params.items() 
                         if isinstance(v, (int, str, float))}
    }

    logging.info(f"Mean: {result['mean']}")
    logging.info(f"Std dev: {result['stddev']}")
    logging.info(f"Avg runtime: {result['avgruntime']}")
    logging.info("\n---------------\n")

    return result


def write_results(results, filename):
    with open(filename, mode="w") as file:
        json.dump(results, file, indent=4)
\end{lstlisting}
\captionof{listing}{The source code for \texttt{util.py}, %
a \Python\ module containing utility functions used for running the simulations.}

\subsection{Source code for the deep Galerkin method}\phantom{h}
\label{sec:source_dgm}

\begin{lstlisting}
import torch

# Computes an approximation of 𝔼[|v(0,X) - φ(X)|²] using the rows of x as independent
# realizations of the random variable X
def dgm_initial_loss(v, φ, x, dev):
    t = torch.zeros(x.shape[0],1).to(dev)
    u = v(torch.cat((t,x), axis=1))
    return (u - φ(x)).square().mean()

# Computes an approximation of 𝔼[|∂v/∂t(T,X) - Δₓv(T,X) - f(v(T,X))|²] using the 
# rows of t and x as independent realizations of the random variables T and X, 
# respectively
def dgm_dynamic_loss(v, ρ, f, t, x, dev):
    d = x.shape[1]
    
    # Split x up into components
    x_comp = [x[:,i] for i in range(d)]
    
    # Compute v(t,x)
    u = v(torch.stack([t] + x_comp, 1))
    
    # Compute ∂v/∂t(t,x)
    u_t = torch.autograd.grad(u, t, torch.ones_like(u), create_graph=True)[0]
    
    # Compute Δₓv(t,x)
    Δu = torch.zeros(x.shape[0]).to(dev)
    for i in range(d):
        # Compute ∂v/∂xᵢ(x,t)
        u_xi = torch.autograd.grad(u, x_comp[i], torch.ones_like(u), 
                                   create_graph=True)[0]
        # Compute ∂²v/∂xᵢ²(x,t)
        u_xixi = torch.autograd.grad(u_xi, x_comp[i], torch.ones_like(u_xi), 
                                     create_graph=True)[0]
        Δu += u_xixi
        
    return (u_t - ρ * Δu - f(u)).square().mean()

# Applies Xavier initialization to the module m (if it is a Linear module)
def init_weights(m):
    if isinstance(m, torch.nn.Linear):
        torch.nn.init.xavier_uniform_(m.weight)
        torch.nn.init.normal_(m.bias, std=0.1)

def train(φ, ρ, f, d, T, model, steps, batch_size, lr, stddev, lr_decay, dev):
    with torch.no_grad():
        # Initialize the weights and biases
        model.apply(init_weights)
        # Set up the Adam optimizer
        optimizer = torch.optim.Adam(model.parameters(), lr=lr)
        # Set up the exponential learning rate decay
        scheduler = torch.optim.lr_scheduler.ExponentialLR(optimizer, lr_decay)

    for step in range(steps):
        with torch.no_grad():
            # Generate the collocation points
            sd = stddev((steps - step) / steps)
            x = (torch.randn(batch_size,d) * sd).to(dev)
            x.requires_grad_()
            t = T * torch.rand(batch_size).to(dev)
            t.requires_grad_()

        optimizer.zero_grad()
        # Compute the loss
        initial_loss = dgm_initial_loss(model, φ, x, dev)
        dynamic_loss = dgm_dynamic_loss(model, ρ, f, t, x, dev)
        loss = initial_loss + dynamic_loss
        # Compute the gradients
        loss.backward()
        # Apply the change to the weights and biases of model
        optimizer.step()
        # Decrease the learning rate
        scheduler.step()
\end{lstlisting}
\captionof{listing}{The source code for \texttt{dgm.py}, %
a \Python\ module containing the core of the implementation of the deep %
Galerkin method.}

\medskip

\begin{lstlisting}
import torch

# Returns a fully connected feed-forward model with the given layer 
# dimensions in which the affine transformation in each hidden layer 
# is followed by the given activation function
def build_model(dev, activation, *layer_dims):
    return lambda: torch.nn.Sequential(
            *[l for (i, j) in zip(layer_dims, layer_dims[1:-1])
                for l in [torch.nn.Linear(i,j), activation()]],
            torch.nn.Linear(layer_dims[-2], layer_dims[-1])
        ).to(dev)
\end{lstlisting}
\captionof{listing}{The source code for \texttt{dgm\_util.py},  %
a \Python\ module containing a utility function for building the DNNs used in  %
the simulations for the deep Galerkin method.}

\medskip

\begin{lstlisting}
import sys
import logging
import torch

from dgm import train
from util import JobConfig, run_jobs, report, write_results
from dgm_util import build_model

# Use the GPU if available
dev = torch.device("cuda" if torch.cuda.is_available() else "cpu")

T = .5
ρ = 1.

# The initial values of the PDE
def φ1(x):
    return (1. + x.square().sum(axis=1, keepdim=True)).sqrt()

def φ2(x):
    return 2. / (4. + x.square().sum(axis=1, keepdim=True))

def φ3(x):
    return torch.atan(x.square().sum(axis=1, keepdim=True).sqrt() / 2.)

# The nonlinearity of the PDE
def f(y):
    return torch.sin(y)

def evaluator(d):
    return lambda model: model(torch.tensor([[T] + [0.0] * d]).to(dev)).item()

configs1 = [
    JobConfig(
        build_model = build_model(dev, torch.nn.Mish, d + 1, d + 40, d + 40, 
                                  d + 40, 1),
        runs = 20,
        train = lambda model, train_params: 
            train(φ1, ρ, f, model=model, dev=dev, **train_params),
        train_params = {
            "d": d,
            "T": T,
            "steps": 350 if d<20 else 
                     500 if d<50 else 
                     750,
            "batch_size": 256, 
            "lr": 0.05 if d<20 else 
                  0.025 if d<50 else 
                  0.01 if d<200 else 
                  0.005,
            "lr_decay": 0.98 if d<20 else 
                        0.99 if d<50 else 
                        0.995,
            "stddev": lambda p: 2**(3. * p -.5)
        },
        evaluate = evaluator(d)
    ) for d in [1, 2, 5, 10, 20, 50, 100, 200]
]

configs2 = [
    JobConfig(
        build_model = build_model(dev, torch.nn.Mish, d + 1, d + 40, d + 40, 
                                  d + 40, 1),
        runs = 20,
        train = lambda model, train_params: 
            train(φ2, ρ, f, model=model, dev=dev, **train_params),
        train_params = {
            "d": d,
            "T": T,
            "steps": 500 if d<10 else 
                     750 if d<50 else 
                     1000,
            "batch_size": 256, 
            "lr":  0.05 if d<10 else 
                   0.025 if d<50 else 
                   0.01 if d<200 else 
                   0.005,
            "lr_decay": 0.99 if d<10 else 
                        0.995 if d<50 else 
                        0.998,
            "stddev": lambda p: 2**(3. * p -.5)
        },
        evaluate = evaluator(d)
    ) for d in [1, 2, 5, 10, 20, 50, 100, 200]
]

configs3 = [
    JobConfig(
        build_model = build_model(dev, torch.nn.Mish, d + 1, d + 40, d + 40, 
                                  d + 40, 1),
        runs = 20,
        train = lambda model, train_params: 
            train(φ2, ρ, f, model=model, dev=dev, **train_params),
        train_params = {
            "d": d,
            "T": T,
            "steps": 350 if d<20 else 
                     500 if d<50 else 
                     750,
            "batch_size": 256, 
            "lr": 0.05 if d<20 else 
                  0.025 if d<50 else 
                  0.01 if d<200 else 
                  0.005,
            "lr_decay": 0.98 if d<20 else 
                        0.99 if d<50 else 
                        0.995,
            "stddev": lambda p: 2**(3. * p -.5)
        },
        evaluate = evaluator(d)
    ) for d in [1, 2, 5, 10, 20, 50, 100, 200]
]

logging.basicConfig(format='%(message)s', level=logging.DEBUG, stream=sys.stdout)

results1 = run_jobs(configs1, lambda: torch.manual_seed(0), report)
write_results(results1, "dgm_1.json")

results2 = run_jobs(configs2, lambda: torch.manual_seed(0), report)
write_results(results2, "dgm_2.json")

results3 = run_jobs(configs3, lambda: torch.manual_seed(0), report)
write_results(results3, "dgm_3.json")
\end{lstlisting}
\captionof{listing}{The source code for \texttt{dgm\_examples.py}, %
a \Python script that runs the simulations for the deep Galerkin method %
and writes the results to files.}

\medskip

\subsection{Source code for the deep splitting method}\phantom{h}
\label{sec:source_ds}

\begin{lstlisting}
import torch
import math

# Computes an approximation of 𝔼[|φ(√(2ρT) W + ξ) + Tf(φ(√(2ρT) W + ξ)) - N(ξ)|²] 
# with W a standard normal random variable using the rows of x as independent 
# realizations of the random variable ξ
def loss(φ, ρ, model, f, T, x, dev):
    with torch.no_grad():
        W = torch.randn_like(x).to(dev)
        z = φ(x + math.sqrt(2 * ρ * T) * W)
        y = z + T * f(z)
    return (y - model(x)).square().mean()

# Applies Xavier initialization to the module m (if it is a Linear module)
def init_weights(m):
    if isinstance(m, torch.nn.Linear):
        torch.nn.init.xavier_uniform_(m.weight)
        torch.nn.init.normal_(m.bias)

def train(φ, ρ, f, d, T, models, lr, lr_decay, steps, batch_size, dev):
    n = len(models)
    for i in range(n):
        # Train model[i]
        with torch.no_grad():
            # Initialize the weights and biases
            models[i].apply(init_weights)
            # Set up the Adam optimizer
            optimizer = torch.optim.Adam(models[i].parameters(), lr=lr)
            # Set up the exponential learning rate decay
            scheduler = torch.optim.lr_scheduler.ExponentialLR(optimizer, lr_decay)
        
        models[i].train()
        for step in range(steps):
            # Generate the collocation points
            with torch.no_grad():
                x = math.sqrt(2 * (n-i-1) * T/n) * torch.randn(batch_size, d).to(dev)
            optimizer.zero_grad()
            # Compute the loss
            loss_value = loss(φ if i==0 else models[i-1], ρ, models[i], f, T/n, x, dev)
            # Compute the gradients
            loss_value.backward()
            # Apply changes to weights and biases of model[i]
            optimizer.step()
            # Decrease the learning rate
            scheduler.step()

        models[i].eval()
\end{lstlisting}
\captionof{listing}{The source code for \texttt{ds.py}, %
a \Python\ module containing the core of the implementation of the deep %
splitting method.}

\medskip

\begin{lstlisting}
import torch

# Returns a list of N fully connected feed-forward models. The first N-1 models have
# the given layer dimensions, in which the affine transformation in each hidden layer 
# is followed by the given activation function. Optionally batch normalization is 
# inserted before each activation function. The N-th model is a single affine linear 
# transformation.
def build_models(dev, N, activation, batch_normalization, *layer_dims):
    def builder():
        stages = [((
            torch.nn.Linear(i,j), 
            *((torch.nn.BatchNorm1d(j),) if batch_normalization else ()), 
            activation()
            ) for (i, j) in zip(layer_dims, layer_dims[1:-1])
        ) for _ in range(N-1)]
        return [torch.nn.Sequential(
            *[m for stage in stages[i] for m in stage],
            torch.nn.Linear(layer_dims[-2], layer_dims[-1]),
            *((torch.nn.BatchNorm1d(layer_dims[-1]),) if batch_normalization else ()), 
        ).to(dev)
        for i in range(N-1)
        ] + [torch.nn.Linear(layer_dims[0],layer_dims[-1]).to(dev)]
    return builder
\end{lstlisting}
\captionof{listing}{The source code for \texttt{ds\_util.py}, %
a \Python\ module containing a utility function for building the DNNs used in %
the simulations for the deep splitting method.}

\begin{lstlisting}
import sys
import logging
import torch

from ds import train
from util import JobConfig, run_jobs, report, write_results
from ds_util import build_models

# Use the GPU if available
dev = torch.device("cuda" if torch.cuda.is_available() else "cpu")

T = .5
ρ = 1.

# The initial values of the PDE
def φ1(x):
    return (1. + x.square().sum(axis=1, keepdim=True)).sqrt()

def φ2(x):
    return 2. / (4. + x.square().sum(axis=1, keepdim=True))

def φ3(x):
    return torch.atan(x.square().sum(axis=1, keepdim=True).sqrt() / 2.)

# The nonlinearity of the PDE
def f(y):
    return torch.sin(y)

def evaluator(d):
    return lambda models: models[-1](torch.zeros(1,d).to(dev)).item()

def train_params(d):
    return {
        "d": d,
        "T": T,
        "steps": 350 if d<20 else 
                 500 if d<100 else 
                 750 if d<200 else 
                 1000 if d<1000 else 
                 1250,
        "batch_size": 256, 
        "lr": 0.2,
        "lr_decay": 0.985 if d<20 else 
                    0.99 if d<100 else 
                    0.995 if d<1000 else 
                    0.998
    }

def trainer(φ):
    return lambda models, train_params: train(φ, ρ, f, models=models, 
                                              dev=dev, **train_params)

configs = [
    [
        JobConfig(
            build_model = 
                build_models(dev, 30, torch.nn.ELU, True, d, d + 50, d + 50, 1),
            runs = 20,
            train = trainer(φ),
            train_params = train_params(d),
            evaluate = evaluator(d)
        ) for d in [1, 2, 5, 10, 20, 50, 100, 200, 500, 1000] 
    ] for φ in [φ1, φ2, φ3]
]

logging.basicConfig(format='%(message)s', level=logging.DEBUG, stream=sys.stdout)

results1 = run_jobs(configs[0], lambda: torch.manual_seed(0), report)
write_results(results1, "ds_1.json")

results2 = run_jobs(configs[1], lambda: torch.manual_seed(0), report)
write_results(results2, "ds_2.json")

results3 = run_jobs(configs[2], lambda: torch.manual_seed(0), report)
write_results(results3, "ds_3.json")
\end{lstlisting}
\captionof{listing}{The source code for \texttt{ds\_examples.py}, %
a \Python\ script that runs the simulations for the deep splitting method %
and writes the results to files.}

\subsubsection*{Acknowledgements}
This work has been funded by the Deutsche Forschungsgemeinschaft 
(DFG, German Research Foundation) 
under Germany's Excellence Strategy EXC 2044 -- 390685587, 
Mathematics M\"unster: Dynamics -- Geometry -- Structure
and through the research grant HU1889/7-1.
This project has been partially supported by the startup fund project
of Shenzhen Research Institute of Big Data under grant No.~T00120220001.

\bibliographystyle{acm}
\bibliography{refnew}

\end{document}